\documentclass{ajour}

\usepackage{latexsym}
\usepackage{amssymb}
\usepackage{amsmath}

\DeclareMathOperator{\+}{\text{ \& }}

\DeclareMathOperator{\finsub}{\text{ is a finite subset of }}

\newcommand{\set}[1]{\left\{\,#1\,\right\}}
\newcommand{\ssim}{\sim \!}
\begin{document}

%%%%% To be entered at Academic Press: =====>>
%
% Uncomment line below only when doing final typesetting,
%\finaltypesetting
% \journame{}
% \articlenumber{}
% \yearofpublication{}
% \volume{}
% \cccline{}
% \received{}
% \revised{}
% \accepted{}

% communication line, use: \commline{Communicated by...}
% \commline{Communicated by... }

\authorrunninghead{Michael Hardy}
\titlerunninghead{Scaled Boolean Algebras}

%% To set particular page number:
%\setcounter{page}{261} %%

%% <<== End of commands to be entered at Academic Press

\title{Scaled Boolean Algebras}

\author{Michael Hardy}
\affil{Department of Mathematics \\
Massachusetts Institute of Technology \\
Cambridge, Massachusetts 02139}

\email{hardy@math.mit.edu}

\abstract{Scaled Boolean algebras are a category of mathematical
objects that arose from attempts to understand why the conventional
rules of probability should hold when probabilities are construed,
not as frequencies or proportions or the like, but rather as degrees
of belief in uncertain propositions.  This paper separates the
study of these objects from that not entirely mathematical problem
that motivated them.  That motivating problem is explicated in the
first section, and the application of scaled Boolean algebras to
it is explained in the last section.  The intermediate sections
deal only with the mathematics.  It is hoped that this isolation
of the mathematics from the motivating problem makes the
mathematics clearer.}

% text should be lower case, unless caps are necessary for meaning

\begin{article}

\section{Introduction}

\subsection{First glimpse}

{\bf Q:} Why should the conventional rules of probability hold
when probabilities are assigned, not to {\em events} that are
{\em random} according to their {\em relative frequencies}
of occurrence, but rather to {\em propositions} that are
{\em uncertain} according to the degree to which they are
supported by the evidence?

{\bf A:} Because probability measures should preserve both the
logical partial ordering of propositions (ordered by logical
implication) and the operation of relative negation.

The explanation and justification of this proposed answer are
not entirely mathematical and appear in \S\ref{bayes} --- the
last section of this paper.
Our main concern will be the mathematical theory that the answer
motivates: the theory of mappings that, like probability measures
and Boolean isomorphisms, preserve partial orderings and certain
kinds of relative complementations.

\subsection{What are scales and what do they measure?}

A scale, as we shall use that word, amounts to a partially
ordered set with what will be called an ``additive relative
complementation.''  In some ways, these behave like
lattice-theoretic relative complementation, although some
of the posets on which they are defined are not lattices,
and when they are lattices, the additive relative
complementation and the lattice-theoretic relative
complementation usually differ.
Additive relative complementation shares the following two
properties with lattice-theoretic relative complementation:
If $\alpha\leq\beta\leq\gamma$, then the additive
complement $\ssim\beta_{[\alpha,\gamma]}$ of $\beta$ relative
to the interval $[\alpha,\gamma]$ goes down from $\gamma$ to
$\alpha$, as $\beta$ goes up from $\alpha$ to $\gamma$, and
the additive complement of the additive complement of $\beta$
is $\beta$ (both complements being taken relative to the same
interval $[\alpha,\gamma]$).  Indeed, when the scale is a
Boolean algebra, then the additive relative complementation
and the lattice-theoretic relative complementation are the
same.  But when the scale is an interval $[\alpha,\gamma]$
on the real line, which is a lattice with no lattice-theoretic
relative complementation, then the additive relative
complementation is
$\beta\mapsto \ \ssim\beta_{[\alpha,\gamma]}=\gamma-\beta+\alpha$.

In lattices (and in particular, in Boolean algebras)
two operations with which conventional lattice-theoretic
relative complementation neatly meshes are each other's duals:
the meet and the join.  In scales, two operations with which
additive relative complementation neatly meshes are again
duals of each other: We shall call them addition ($+$) and
dual-addition ($\oplus$).  When the scale is a Boolean
algebra, then addition and dual addition coincide exactly
with meet and join.  When the scale is the interval
$[0,1]$ on the real line, then addition is ordinary addition
restricted to pairs of numbers whose ordinary sum is
within the interval $[0,1]$, and the dual addition is
$(\alpha,\beta)\mapsto\alpha+\beta-1$, again restricted
to pairs of numbers for which the value of that operation
is within the interval.
Two de-Morganesque laws relate additive relative
complementation to addition and dual addition,
and a ``modular law'' says
$(\zeta+\eta)\oplus\theta=\zeta+(\eta\oplus\theta)$
when $\zeta,\eta,\theta$ are suitably situated.  This modular
law is reminiscent of the law defining modular lattices,
which says $(x\vee y)\wedge z=x\vee(y\wedge z)$ whenever
$x\leq z$ (where ``$\wedge$'' and ``$\vee$'' are meet and join).

All scales live somewhere between (i.e., {\em inclusively} between)
the Boolean case and various sorts of linearly ordered cases.

Probability, construed somewhat liberally, will be
measured on linearly ordered scales --- we will allow some
generalizations of probability measures to take values in
other linearly ordered scales than the interval $[0,1]$ on
the real line.  These generalizations will, like probability
measures, preserve the partial order and the relative
complementation.  I think I chose the word ``scale''
because of talk among followers of Edwin Jaynes about
``scales'' on which probability can be measured.
They were inspired by a theorem from the physicist
Richard T.~Cox's book {\em Algebra of Probable Inference}
\cite{cox} that under certain semi-plausible (and too strong,
in my view) assumptions about the behavior of probabilities,
construed as degrees of belief in uncertain propositions, all
such scales must be isomorphic to the unit interval
$[0,1]\subseteq\mathbb{R}$ with its usual furniture -- addition,
multiplication, and linear ordering -- and probability measures
must conform to the usual addition and multiplication rules.

This paper resulted from my attempt to understand Cox's
arguments, but what we do here will be quite different
from those, to say the least.

\subsection{Frequentism and Bayesianism}

``Frequentists'' assign probabilities to random events
according to their relative frequencies of occurrence,
or to subsets of populations as proportions of the whole.
``Bayesians'', on the other hand, assign probabilities
to uncertain propositions according to the degree to which
the evidence supports them.  Frequentists treat
probabilities as {\em intrinsic} to the object of study,
and Bayesians treat them as {\em epistemic}, i.e., conditional
on one's {\em knowledge} of the object of study.
Frequentist and Bayesian methods of statistical inference
differ, and their relative merits have been debated for decades.

Here is a poignant example of a problem whose space
of feasible solutions changes when the Bayesian outlook
replaces the frequentist one.
The respective prices of three kinds of gadgets are \$20,
\$21, and \$23.  Records of the gross receipts of Acme Gadgets
for the year 2099 show that customers bought 3,000,000,000,000,000
gadgets during that year, spending \$66,000,000,000,000,000
on them, so that they spent an average of \$22 per gadget.
One of those gadgets sits in an unopened box on your desk.
It is just as likely to be any of the 3,000,000,000,000,000
gadgets as it is to be any of the others. \\
{\bf Ambiguous question:} Given this information, what are
the probabilities that this gadget is of the first, second,
and third kinds? \\
{\bf A frequentist way to construe the question:}
What proportions $p_1$, $p_2$, $p_3$, of the first, second,
and third kinds, respectively, were purchased? \\
{\bf A Bayesian way to construe the question:}
What degrees of belief $p_1, p_2, p_3$ should we assign
to the propositions that {\em this particular} gadget is
respectively of the first, second, or third kind, if we
have only the information given above and no more?

Under the frequentist construction of the question,
the feasible solutions are the convex combinations
of the two extreme solutions
\begin{eqnarray*}
(p_1,p_2,p_3) & = & (0, \ 1/2, \ 1/2), \\
(p_1,p_2,p_3) & = & (1/3, \ 0, \ 2/3).
\end{eqnarray*}
But both of these extreme solutions are incorrect under
the Bayesian construction of the question!  The first
extreme solution says $p_1=0$.  By the Bayesian construction,
this would mean that, on the meager evidence given, we can be
sure that the gadget in the box is not of the first kind.
But obviously we cannot.  The other extreme solution says
$p_2=0$, and is defeated by the same consideration.
Whether there is any solution, and whether there is only
one solution, on the Bayesian construction, is a subtler
question.  (Edwin Jaynes' ``principle of maximum entropy,''
proposed in \cite{jaynes}, entails that the rationally
justified degrees of belief are those that maximize the
entropy $\sum_{i=1}^3 -p_i\log p_i$ subject to the
constraints $p_1, p_2, p_3 \geq 0$, \ \ $p_1+p_2+p_3=1$,
and $\$20p_1+\$21p_2+\$23p_3=\$22$.)

\subsection{Why the conventional rules of probability?}

The conventional mathematical rules of probability include
additivity and definitions and other characterizations of
conditional probability.  To ``pure'' mathematicians, these
are merely axioms or definitions.  To frequentists, finite
additivity and the definition of conditioning on events of
positive probability are trivially true propositions about
frequencies or proportions.  But to Bayesians, the rules of
probability are problematic.

In \cite{cox} Cox addressed the question of why finite additivity
-- the ``sum rule'' -- and the conventional definition of
conditional probability -- the ``product rule'' -- should be
adhered to if probabilities are taken to be degrees of belief
rather than proportions or frequencies or the like.

\subsection{How and why this work came about}

I set out to recast Cox's argument in a more abstract form.
However, the project went in a direction of its own choosing.
Ultimately I had a different, but not more abstract, argument
(call it the ``concrete version'') for a similar but more
extensive set of conclusions, relying on a different
(and more plausible, in my view) set of assumptions about
rational degrees of belief.
Then I set out to recast {\em that} argument in more abstract
language (the ``more abstract version'') in order to separate
the part of it that is purely mathematical from the rest.
To that end, I conceived the idea of a scale.

The concrete version of the argument rests on an explanation
of why probabilities, construed as degrees of belief in uncertain
propositions, should be assigned in a way that preserves the
logical partial ordering and relative complementation of
propositions.  That explanation is not entirely mathematical,
and appears in \S\ref{bayes}.

\subsection{Relation of this to earlier work}\label{savages}

The concrete version overlaps with some earlier work of
Leonard Jimmie Savage \cite{savage}, Terrence Fine \cite{fine},
Peter C.~Fishburn \cite{fishburn}, Bruno de Finetti \cite{definetti},
Bernard O.~Koopman \cite{koopman}, and Kraft, Pratt, and Seidenberg
\cite{kraftetal} that I had initially largely ignored because
those authors seemed to be assuming a weak sort of additivity
as an axiom rather than trying to prove that probabilities should
be assigned additively.  That weak additivity statement appears
here as Lemma~\ref{basiclemma}.  I found that I had rediscovered
the result that appears here as Theorem~\ref{additivity}.

However, a seemingly trivial change in emphasis makes it possible
to go considerably beyond where those authors left off, and hence
to define the concept of a scale and its addition, dual addition,
and additive relative complementation.
Those authors considered two orderings of propositions:
the first ``$\leq$'' is the usual logical partial ordering, so that
$x\leq y$ if $x$ logically entails $y$.  The second, a linear
ordering, ``$\preceq$'', is a comparative probability ordering,
so that ``$x\preceq y$'' means $x$ is no more probable than $y$.
Those authors assumed:
\begin{equation}\label{prob-ordering}
\begin{split}
\text{If }x\leq y\text{ then }x\preceq y,
\text{ so ``}\preceq
\text{'' is a linear extension of ``}\leq\text{''}. \\
\text{If }x\preceq y\text{ and }y\wedge z=0,
\text{ then }x\vee z\preceq y\vee z.\quad
\text{(weak additivity)}
\end{split}
\end{equation}
The seemingly trivial change in emphasis is from an ordering
of propositions by probabilities to a suitably well-behaved
mapping, called a scaling, from a Boolean algebra of propositions
into a partially ordered space of generalized ``probabilities.''

Kraft, Pratt, and Seidenberg \cite{kraftetal} found that a
condition called ``strong additivity'' is sufficient to
guarantee that a comparative probability ordering $\preceq$
has an ``agreeing measure,'' i.e., a probability measure (in the
usual sense of the term) $\mu$, such that $x\preceq y$ holds if
and only if $\mu(x)\leq\mu(y)$.  We give an enormously simpler
sufficient condition called ``divisibility'' in \S\ref{divisibility}.
In \cite{scott}, Dana Scott proved Kraft, Pratt and Seidenberg's
result by a more generally applicable method.  Scott showed that
his method can be applied not only to probabilities, but also to
other things naturally measured on partially ordered scales.

Unlike the authors cited above, we also consider infinite
Boolean algebras.  When the Boolean algebra that is the domain
of a scaling is infinite, it makes sense to speak of continuity
or discontinuity of a scaling at a particular Boolean homomorphism.
We shall see that continuity at all homomorphisms whose
kernels are principal ideals is the same as complete additivity,
and continuity at all 2-valued homomorphisms entails a kind
of Archimedeanism.

\newpage
\section{Boolean algebras and scales}\label{BA}

\subsection{Boolean algebras}

\begin{definition}\label{dfnboole}
A {\bf Boolean algebra} consists of an underlying
set $\mathbb{A}$ with two distinguished elements
$0 \neq 1$, two binary operations
$(x,y) \mapsto x \wedge y=$ the ``meet'' of $x$ and $y$, and
$(x,y) \mapsto x \vee y=$ the ``join'' of $x$ and $y$, and
a unary operation $x \mapsto \ssim x=$ the ``complement'' of $x$,
satisfying the following algebraic laws (which are the same laws
that are obeyed by the logical connectives ``and'', ``or'',
``not'' or the operations of intersection, union, and
complementation of sets):
For $x,y,z\in\mathbb{A}$ we have
\begin{equation}\label{booleIDs}
\begin{array}{cc}
x \wedge y=y \wedge x   &   x \vee y=y \vee x \\
(x \wedge y) \wedge z=x \wedge (y \wedge z) &
       (x \vee y) \vee z=x \vee (y \vee z) \\
x \wedge (y \vee z)=(x \wedge y) \vee (x \wedge z) \quad &
\quad x \vee (y \wedge z)=(x \vee y) \wedge (x \vee z) \\
x \wedge x=x & x \vee x=x \\
 \multicolumn{2}{c}{\sim \ssim x=x} \\
\ssim(x \wedge y)=(\ssim x) \vee (\ssim y) &
  \ssim(x \vee y)=(\ssim x) \wedge (\ssim y) \\
x \wedge \ssim x=0 & x \vee \ssim x=1 \\
x \wedge 1=x & x \vee 0=x \\
x \wedge 0=0 & x \vee 1=1 \\
\ssim 1=0 & \ssim 0=1
\end{array}
\end{equation}
\end{definition}
By a ``convenient abuse of language'' we shall
use the same symbol $\mathbb{A}$ to refer either
to a Boolean algebra or to its underlying set.

Every Boolean algebra has a natural partial ordering ``$\leq$'' 
defined by $x \leq y$ iff $x \wedge y=x$, or equivalently
$x \vee y=y$.  With this ordering, $x \wedge y$ and $x \vee y$
are respectively the infimum and the supremum of the set $\set{x,y}$.
The largest and smallest elements of $\mathbb{A}$
are respectively 1 and 0.

\begin{definition}
For $a,b,\in\mathbb{A}$, if $a \leq b$ then
the {\bf complement} of any
$x\in[a,b]=\set{u : a \leq u \leq b}$ {\bf relative}
to the interval $[a,b]$ is
$$
\ssim x_{[a,b]}=a \vee (b \wedge \ssim x)
=b \wedge (a \vee \ssim x).
$$
\end{definition}
\begin{proposition}
For $a,b\in\mathbb{A}$, if $a \leq b$ then the interval
$[a,b]=\set{x : a \leq x \leq b}$
is a Boolean algebra whose meet and join operations are
the restrictions to $[a,b]$ of the meet and join operations
of $\mathbb{A}$, and whose complementation is relative
to this interval.
\end{proposition}
The proof is a quick exercise.

When $[a,b]=[0,1]=\mathbb{A}$ then relative complementation
coincides with ordinary complementation.  An interval
$[a,b]\subsetneq\mathbb{A}$ with this structure is not a
``subalgebra'' of $\mathbb{A}$ because its complementation
differs from that of $\mathbb{A}$.  An interval $[a,b]$ with
this structure will be called a ``relative Boolean algebra.''

If $[a,b] \subseteq [c,d] \subseteq [0,1]=\mathbb{A}$ then
the intervals $[c,d]$ and $[0,1]$ are both larger Boolean
algebras of which $[a,b]$ is a subinterval.  Should one take
\mbox{$\ssim x_{[a,b]}$} to be $a \vee (b \wedge \ssim x)$
or \mbox{$a \vee (b \wedge \ssim x_{[c,d]})$}?
A straightforward computation shows that either 
yields the same result.

The algebraic laws defining a Boolean algebra are those
obeyed by the logical connectives ``and'', ``or'', ``not''
that connect propositions.  The ``$0$'' and ``$1$'' in
a Boolean algebra correspond respectively to propositions
known to be false, and propositions known to be true.
The relation $x \leq y$ corresponds to the
statement that it is known that if $x$ is true then so is
$y$, although the truth values of $x$ and $y$ may be uncertain.
If $0\leq x<y\leq 1$ then $y$ is closer to being known
to be true than $x$ is; $x$ is closer to being known to be
false than $y$ is.  We shall argue that the definition of
``scaling'' that will follow, captures laws that should
be obeyed by any rational assignment of degrees of belief
to propositions.

\subsection{Basic scalings and scales}

A Boolean algebra has a partial ordering and a relative
complementation.  A ``scale'' is a more general kind of
object with a partial ordering and an ``additive'' relative
complementation (although, as we shall see, the latter fails
to be everywhere-defined in some cases).  One of the simplest
examples of a scale that is not a Boolean algebra is the
interval $[0,1]\subseteq\mathbb{R}$, in which the additive
complement of $\beta\in [\alpha,\gamma]$ relative to the
interval $[\alpha,\gamma]$ is
\mbox{$\ssim\beta_{[\alpha,\gamma]}
=\gamma-\beta+\alpha\in[\alpha,\gamma]$}.
A ``scaling'' is a mapping from one scale to another that
preserves the partial ordering and the additive relative
complementation.  A ``basic'' scaling is one whose domain
is a Boolean algebra.  Among the simplest basic scalings are
Boolean isomorphisms and finitely additive measures.

\begin{definition}\label{BasicDefinitions}
\begin{enumerate}
\item\label{BasicScalingDfn}
A {\bf basic scaling}
is a mapping $\rho$ from a Boolean algebra $\mathbb{A}$ into
any partially ordered set, that (a) is strictly increasing,
and (b) preserves relative complementations.
These two conditions mean:
 \begin{enumerate}
  \item\label{firstcriterion}
    For $x,y \in \mathbb{A}$, if $x<y$ then $\rho(x)<\rho(y)$, and
  \item\label{secondcriterion}
    For $x,y \in [a,b] \subseteq \mathbb{A}$, if $\rho(x)<\rho(y)$
    then $\rho(\ssim x_{[a,b]})>\rho(\ssim y_{[a,b]})$,
    and if $\rho(x)=\rho(y)$ then
    $\rho(\ssim x_{[a,b]})=\rho(\ssim y_{[a,b]})$.
    (In particular, if $a=0$, $b=1$, and $\rho(x)<\rho(y)$,
    then $\rho(\ssim x)>\rho(\ssim y)$ and similarly if ``$=$''
    replaces ``$<$''.)
 \end{enumerate}
(A ``scaling'' is either a basic scaling as defined here,
or a more general kind of scaling to be defined in
\S\ref{Difficulty}.)
\item
A {\bf scaled Boolean algebra} is a Boolean algebra endowed
with a basic scaling.  If $\rho$ is a basic scaling on a Boolean
algebra $(\mathbb{A},0,1,\wedge,\vee,\ssim)$ then
$(\mathbb{A},0,1,\wedge,\vee,\ssim,\rho)$ is a scaled
Boolean algebra.
\item\label{ScaleDfn}
A {\bf scale} is the image
$\mathcal{R}=\set{\rho(x):x \in \mathbb{A}}$ of a basic scaling
$\rho$, with its partial ordering and an additive relative
complementation that it inherits from the scaling
$\rho:\mathbb{A}\rightarrow\mathcal{R}$.  The precise
definition of this inherited additive relative complementation
involves some perhaps unexpected complications, and appears in
\S\ref{Difficulty}.  Lower-case Greek letters will usually be
used for members of $\mathcal{R}$, except that the minimum and
maximum members of $\mathcal{R}$ will be called 0 and 1
respectively.
\end{enumerate}
\end{definition}
Clearly the restriction of a basic scaling to a relative
Boolean algebra $[a,b]$ is also a basic scaling.

\subsection{The difficulty with relative complementation}\label{Difficulty}

By Definition~\ref{BasicDefinitions} (1ii), a scaling
induces a sort of complementation on a scale -- one may
unambiguously define $\ssim\rho(x)$ to be $\rho(\ssim x)$.
This operation is strictly decreasing and is its own inverse:
For $\alpha,\beta$ in a scale $\mathcal{R}$, we have
$\ssim\alpha>\ssim\beta$ if $\alpha<\beta$, and
\mbox{$\sim\ssim\alpha=\alpha$}.
It may be tempting to think it is just as easy to define
an induced complementation {\em relative} to an interval.
Here is the difficulty.  Suppose $x,y \in [a,b]\subseteq\mathbb{A}$.
Although the definition of ``basic scaling'' says that
if $\rho(x)=\rho(y)$ then
\mbox{$\rho(\ssim x_{[a,b]})=\rho(\ssim y_{[a,b]})$},
the extension to relative complements requires
something stronger.  We need to know the following fact.

{\em Suppose $\rho$ is a basic scaling.
If $x\in[a,b]$ and $y\in[c,d]$ and
$$
\rho(a)=\rho(c) \ \ \leq \ \ \rho(x)=\rho(y)
\ \ \leq \ \ \rho(b)=\rho(d)
$$
then $\rho(\ssim x_{[a,b]})=\rho(\ssim y_{[c,d]})$.}

\noindent
More economically stated: $\rho(\ssim x_{[a,b]})$
depends on $x$, $a$, and $b$ only through $\rho(x)$, $\rho(a)$,
and $\rho(b)$.  The proof appears in \S\ref{relcomp}.  
This proposition makes the following definition unambiguous.
\begin{definition}
Additive relative complementation on $\mathcal{R}$ is given by
$$\ssim\rho(x)_{[\rho(a),\rho(b)]}=\rho(\ssim x_{[a,b]}).$$
\end{definition}
The word ``additive'' is used because of the relationship
between this relative complementation and the operations of
addition and dual addition.  Those operations are introduced
in \S\ref{addition}.
Additive relative complementation may fail to be everywhere-defined,
since it can happen that $\zeta<\eta<\theta$ even while no
$x<y<z$ exist in $\mathbb{A}$ whose respective images under
$\rho$ are $\zeta,\eta,\theta$.  Concrete instances will appear in
\S\ref{examples}.  Scalings for which this happens are not
``divided.''  That concept is defined in \S\ref{divisibility},
which section also prescribes the remedy to this pathology.

\subsection{General definition of scaling}

Here is the definition of ``scaling'' that is more general
than that of ``basic scaling.''
\begin{definition}\label{GeneralScalingDfn}
A {\bf scaling} is a mapping $\rho$ from a scale $\mathcal{R}$
into a partially ordered set, that (a) is strictly increasing,
and (b) preserves additive relative complementations.
These two conditions mean:
 \begin{enumerate}
  \item
    For $\zeta,\eta\in\mathcal{R}$,
    if $\zeta<\eta$ then $\rho(\zeta)<\rho(\eta)$, and
  \item
    For $\zeta,\eta\in [\alpha,\beta]\subseteq\mathcal{R}$,
    if $\rho(\zeta)<\rho(\eta)$ then
    $\rho(\ssim\zeta_{[\alpha,\beta]})
    >\rho(\ssim\eta_{[\alpha,\beta]})$,
    and if $\rho(\zeta)=\rho(\eta)$ then
    $\rho(\ssim\zeta_{[\alpha,\beta]})
    =\rho(\ssim\eta_{[\alpha,\beta]})$.
 \end{enumerate}
\end{definition}
Although in this definition $\rho$ need not be a basic
scaling, its image is nonetheless the image of a basic scaling
$\rho\circ\sigma$, where $\mathcal{R}$ is the image of a
basic scaling $\sigma$ on some Boolean algebra $\mathbb{A}$.
Therefore all images of scalings are images of basic scalings,
and we need not change the definition of ``scale'' as the image
of a basic scaling.
\begin{definition}\label{extensiondefinition}
If $\rho:\mathbb{A}\rightarrow\mathcal{R}$ is a basic scaling
and $\sigma:\mathcal{R}\rightarrow\mathcal{S}$ is a scaling, then
the basic scaling $\sigma\circ\rho:\mathbb{A}\rightarrow\mathcal{R}$
is an {\bf extension} of the basic scaling $\rho$.
\end{definition}
Why an ``extension''?  Because $\sigma$ may map two
incomparable members $\alpha,\beta$ of $\mathcal{R}$
to members $\sigma(\alpha),\sigma(\beta)$ that are comparable
either because they are equal (so that $\sigma$ is not one-to-one)
or because a strict inequality holds between them.  In other
words, $\sigma$ extends the partial ordering by making comparable,
and possibly even equal, things that were incomparable before
the extension.  (Note that the definition precludes unequal
but comparable elements of $\mathcal{R}$ having the
same image under $\sigma$.)

\subsection{Measures on Boolean algebras}

\begin{definition}\label{measures.on.BAs}
A {\bf finitely additive measure} on a Boolean algebra
$\mathbb{A}$ is a mapping
\mbox{$\rho:\mathbb{A}\rightarrow[0,\infty)\subseteq\mathbb{R}$},
satisfying
  \begin{enumerate}
     \item
         For all $x \in \mathbb{A}$, if $x>0$ then $\rho(x)>0$; and
     \item
         For all $x,y \in \mathbb{A}$, if $x \wedge y=0$,
         then $\rho(x \vee y)=\rho(x)+\rho(y)$.
  \end{enumerate}
\end{definition}
The requirement that if $x>0$ then $\rho(x)>0$ excludes analogs
of non-empty sets of measure zero.  We deal with things like
Lebesgue measure by regarding sets that differ only by a set
of measure zero as equivalent to each other, and considering the
Boolean algebra of equivalence classes.

Definition~\ref{measures.too} will generalize
Definition~\ref{measures.on.BAs} by defining the concept of a
finitely additive measure whose domain is an arbitrary scale.

\section{Examples}\label{examples}
\begin{example}\label{isomorphism}
Every isomorphism from one Boolean algebra into another is
a basic scaling; hence every Boolean algebra is a scale.
\end{example}
\begin{example}\label{measurescale}
Every finitely additive measure on a Boolean algebra
is a basic scaling.  Since for any $a\in\mathbb{R}$,
$a>0$, there exist measures whose image is the whole
interval $[0,a]\subseteq\mathbb{R}$, that interval is
a scale with relative complementation given by
$\ssim\beta_{[\alpha,\gamma]}=\gamma-\beta+\alpha$
for $\beta\in[\alpha,\gamma]$.
\end{example}
\begin{example}\label{convexscale}
Let $\mathcal{M}$ be a nonempty convex set of finitely additive
measures on a Boolean algebra $\mathbb{A}$.  Call $x,y\in\mathbb{A}$
equivalent if $\mathcal{M}$ does not separate $x$ from $y$,
i.e., $\mu(x)=\mu(y)$ for every $\mu\in\mathcal{M}$.
Let $\rho(x)$ be the equivalence class to which $x$ belongs.
Say that $\rho(x)<\rho(y)$ if for {\em every} $\mu\in\mathcal{M}$
we have $\mu(x)\leq\mu(y)$ and for {\em some}
$\mu\in\mathcal{M}$ the inequality is strict.
Plainly this is an antisymmetric relation; to show that
it is a strict partial ordering we need to show that it
is transitive.  If $\rho(x)<\rho(y)$ and
$\rho(y)<\rho(z)$ then there exist $\mu,\nu\in\mathcal{M}$
such that $\mu(x)<\mu(y)\leq\mu(z)$ and
$\nu(x)\leq\nu(y)<\nu(z)$.  Convexity implies
$\pi=(\mu+\nu)/2\in\mathcal{M}$, and then we have $\pi(x)<\pi(z)$.
The reader can check that $\rho$ is a scaling.
If $\mathcal{M}\neq\varnothing$ is the set of {\em all} measures
on $\mathbb{A}$ then $\rho$ is just the canonical automorphism of
$\mathbb{A}$.  At the opposite extreme, $\mathcal{M}$ could contain
just one measure and $\rho$ would be a linearly ordered scale.
\end{example}
\begin{example}\label{nonarch}
This example is (1) a simple ``non-Archimedean'' scale;
(2) a scale that is not a lattice; (3) a scaling that is
not countably additive, and (4) a scaling that is
discontinuous at some 2-valued Boolean homomorphisms on its
domain.  Precise definitions of terms needed to understand
these claims will appear in the sequel.

For $A,B\subseteq\mathbb{N}=\set{1,2,3,\dots}$, let
$\left| A \right|\in\set{0,1,2,3,\dots,\aleph_0}$
be the cardinality of $A$ and let
$A\diagdown B=\set{x : x \in A \+ x \not\in B}$.
Call two sets $A,B$ equivalent if
$$
\left| A\diagdown B \right|=\left| B\diagdown A \right|<\aleph_0,
$$
i.e., $A$ can be changed into $B$ by deleting finitely many
members and replacing them by exactly the same number of others.
Let $\rho(A)$ be the equivalence class
to which $A$ belongs.  Say that $\rho(A)<\rho(B)$ if
$\left|A\diagdown B\right|<\left|B\diagdown A\right|$.
Note that $\left|A\diagdown B\right|=\left|B\diagdown A\right|$
does not imply $\rho(A)=\rho(B)$ unless the common value of these
two cardinalities is finite.  If
$\left|A\diagdown B\right|=\left|B\diagdown A\right|=\aleph_0$
then $\rho(A)$ and $\rho(B)$ are not comparable.

What does the poset
$\mathcal{R}=\set{\rho(A) : A \subseteq \mathbb{N}}$ look like?
Every member $\alpha$ of $\mathcal{R}$ except $\rho(\varnothing)$ has
an immediate predecessor, a largest member of $\mathcal{R}$ that is
$<\alpha$.  Call it $\alpha-1$.  Similarly, every member of $\mathcal{R}$
except $\rho(\mathbb{N})$ has an immediate successor $\alpha+1$,
the smallest member of $\mathcal{R}$ that is $>\alpha$.  The range
$\mathcal{R}$ is partitioned into ``galaxies''
$$\set{\dots,\alpha-2, \alpha-1, \alpha, \alpha+1, \alpha+2, \dots},$$
plus an ``initial galaxy''
$$\set{\rho(\varnothing),\rho(\varnothing)+1,\rho(\varnothing)+2,\dots}$$
of ``finite elements'' and a ``final galaxy''
$$\set{\dots,\rho(\mathbb{N})-2,\rho(\mathbb{N})-1,\rho(\mathbb{N})}$$
of ``cofinite elements.''  For any two galaxies that are comparable,
in the sense that any member of one is comparable to any member of
the other, uncountably many other galaxies are between them, and
infinite antichains of galaxies are between them.
(An antichain in a partially ordered set is a set of pairwise
incomparable elements.)

This mapping $\rho$ is a scaling.
For any finite element $\alpha>0$,
one can write $\mathbb{N}$ as a union of subsets whose
images under $\rho$ are $\leq\alpha$.
But $\mathbb{N}$ cannot be written as a union of finitely
many such sets, so we say that $\mathcal{R}$ is a
``non-Archimedean'' scale.

Via this example it is easy to see why a scale has more structure
than its partial ordering.  Single out any typical galaxy, and define
a mapping by $\alpha\mapsto\alpha+1$ if $\alpha$ is in that galaxy,
and $\alpha\mapsto\alpha$ otherwise.  This is clearly a
poset-automorphism, but it is not a scale-automorphism since it
fails to preserve relative complementation.

This scale is not a lattice, i.e., a set of two members of
$\mathcal{R}$ need not have an infimum or a supremum.
If $\alpha\leq$ each of two incomparable members of $\mathcal{R}$,
then so is $\alpha+1$, and similarly for ``$\geq$'' and $\alpha-1$.

In this scale the relative complementation is everywhere-defined.
\begin{remark}
Members of this scale can serve as dimensions of
subspaces of infinite-dimensional separable Hilbert space.
\end{remark}
\end{example}
\begin{example}(A linearly ordered non-Archimedean scale)
\label{fullnonarch}

First we tersely describe this example in language
comprehensible to those who know nonstandard analysis.
Then for others we include a two-page crash course
in that subject.

Let $n$ be an infinite integer.
The \mbox{$\ast$-finitely} additive measure on the
set of all {\em internal} subsets of $\set{1,\dots,n}$
that assigns $1$ to every one-element set is a scaling.
As in the last example, there is one member, $1$, of the
range of this scale such that $\set{1,\dots,n}$ can be
written as the union of subsets each of which is mapped
to something $\leq 1$, but $\set{1,\dots,n}$ cannot
be written as the union of finitely many such sets.
As in the previous example, the scale is partitioned
into uncountably many ``galaxies,''
$\set{\dots,k-1,k,k+1,\dots}$, plus and initial galaxy
$\set{0,1,2,\dots}$ and a final galaxy $\set{\dots,n-2,n-1,n}$.
But in this case any two galaxies are comparable.

Now for the two-page crash course.
(This can be skipped by anyone who wants to take the
previous paragraph on faith.)
The ordered field $\mathbb{R}^\ast$ of nonstandard real
numbers properly includes the real field $\mathbb{R}$.
Like all ordered fields that properly include $\mathbb{R}$,
this field is ``non-Archimedean.''  This term, when applied
to ordered fields, has a simpler definition than that used
by people who study fields of $p$-adic numbers.  It means
that some members $x\neq 0$ of $\mathbb{R}^\ast$ are
infinitesimal, i.e.,
$$\underbrace{\left|x\right|+\cdots+\left|x\right|}_{n \text{ terms}}<1
\text{ for every finite cardinal number } n.$$
The only infinitesimal in $\mathbb{R}$ is 0.
Some other members of $\mathbb{R}^*$ -- the reciprocals
$y$ of the nonzero infinitesimals -- are infinite, i.e.,
$$\underbrace{1+\cdots+1}_{n\text{ terms}}<\left|y\right|
\text{ for every finite cardinal number } n.
$$

The underlying set of the field $\mathbb{R}^*$ is the
image of $\mathbb{R}$ under a mapping $A\mapsto A^*$
from subsets $A$ of $\mathbb{R}$ to subsets of
$\mathbb{R}^*$.  In every case $A\subseteq A^*$, with
equality if and only if $A$ is finite.  Sets of the
form $A^*$ for some $A\subseteq\mathbb{R}$ are called
``standard'' subsets of $\mathbb{R}^*$.  The standard
sets belong to a much larger class of subsets of
$\mathbb{R}^*$ called ``internal'' sets.  Similarly
each function $f:A\rightarrow\mathbb{R}$ extends
to a function $f^*:A^*\rightarrow\mathbb{R}^*$;
these are called ``standard'' functions, and belong
to the much larger class of ``internal'' functions.
Sets and functions that are not internal are ``external.''
Although space limitations forbid defining these precisely
here, their role and importance will become evident from
the following proposition and its accompanying examples.
\begin{proposition}\label{transfer}(The ``transfer principle'')
\begin{enumerate}
\item\label{elementary}
Suppose a proposition that is true of $\mathbb{R}$ can be
expressed via functions of finitely many variables
(e.g.~$(x,y)\mapsto x+y$), relations among finitely
many variables (e.g.~$x\leq y$), finitary logical connectives
such as ``and'', ``or'', ``not'', ``if\dots then \dots'',
and the quantifiers $\forall x\in\mathbb{R}$ and
$\exists x\in\mathbb{R}$.
(For example, one such proposition is
$\forall x\in\mathbb{R} \ \exists y\in\mathbb{R} \ x+y=0.$)
Such a proposition is true in $\mathbb{R}$
if and only if it is true in $\mathbb{R}^*$
when the quantifier $\forall x\in\mathbb{R}^*$
replaces $\forall x\in\mathbb{R}$, and similarly
for ``$\exists$''.
\item\label{starmap}
Suppose a proposition otherwise expressible as simply
as those considered in part (\ref{elementary}) above
mentions some particular sets $A\subseteq\mathbb{R}$.
Such a proposition is true in $\mathbb{R}$ if and only
if it is true in $\mathbb{R}^*$ with each such ``$A$''
replaced by the corresponding ``$A^*$''.
(Here are two examples: (1) The set
$[0,1]^\ast=\set{x\in\mathbb{R}:0\leq x\leq 1}^\ast$ must
be $\set{x\in\mathbb{R}^\ast:0\leq x\leq 1}$, including
not only members of $\mathbb{R}$ between 0 and 1
inclusive, but also members of $\mathbb{R}^*$ that differ
from those by infinitesimals.  To see this, observe that
the sentence
$$\forall x\in\mathbb{R} \
(x\in [0,1] \text{ if and only if } 0\leq x \leq 1)$$
is true in $\mathbb{R}$, and apply the transfer principle.
(2) The set $\mathbb{N}^\ast$ must be a set that has
no upper bound in $\mathbb{R}^\ast$ (since the sentence
expressing the non-existence of an upper bound of $\mathbb{N}$
in $\mathbb{R}$ is simple enough for the transfer principle
to apply to it) and must contain $n+1$ if
it contains $n$, but must not contain anything between $n$
and $n+1$.  Members of $\mathbb{N}^\ast\diagdown\mathbb{N}$
are ``infinite integers.'')
\item\label{internal}
Suppose a proposition otherwise expressible as simply
as those considered in parts (\ref{elementary}) and
(\ref{starmap}) above contains the quantifier
``$\forall A\subseteq\mathbb{R}\dots$'' or
``$\exists A\subseteq\mathbb{R}\dots$''.
Such a proposition is true in $\mathbb{R}$ if and only if
it is true in $\mathbb{R}^*$ after the changes specified
above and the replacement of the quantifiers with
``$[\forall \text{ \underline{internal} }
A\subseteq\mathbb{R}^\ast\dots]$''
and ``$[\exists \text{ \underline{internal} }
A\subseteq\mathbb{R}^\ast\dots]$''.
(Here are three examples: (1) Every nonempty {\em internal} subset
of $\mathbb{R}^\ast$ that has an upper bound in $\mathbb{R}^\ast$
has a least upper bound in $\mathbb{R}^\ast$.  Consequently
the set of all infinitesimals is external.
(2) The well-ordering principle implies every nonempty
{\em internal} subset of $\mathbb{N}^\ast$ has a
smallest member.  Consequently the set
$\mathbb{N}^\ast\diagdown\mathbb{N}$ of all infinite
integers is external.  (3) If $n$ is an infinite integer,
then the set $\set{1,\dots,n}$ (which is not standard!) must be internal.
To prove this, first observe that the following is trivially true:
$$\forall n\in\mathbb{N} \ \exists A\subseteq\mathbb{N} \
\forall x\in\mathbb{N} \ [x\in A \text{ iff } x \leq n].$$
Consequently
$$\forall n\in\mathbb{N}^\ast \ \exists \text{ \underline{internal} }
A\subseteq\mathbb{N}^\ast \ \forall x\in\mathbb{N}^\ast
\ [x\in A \text{ iff } x\leq n].)$$
\item\label{StarFiniteExample}
As with internal sets, so with internal functions:
Replace ``$[\forall f:A\rightarrow\mathbb{R}\dots]$''
with ``$[\forall\text{ \underline{internal} }
f:A^\ast\rightarrow\mathbb{R}^\ast\dots]$'', and
similarly with ``$\exists$'' in place of ``$\forall$''.
(For example: If $n$ is an infinite integer, then the
complement of the image of any internal one-to-one
function $f$ from the infinite set $\set{1,\dots,n}$
into $\set{1,\dots,n}\cup\set{n+1,n+2,n+3}$ has exactly
three members.  Because of the infiniteness of the domain,
the complements of the images of one-to-one functions
from the former set to the latter come in many sizes,
but most of these functions are external.)
\end{enumerate}
\end{proposition}
The last example described in Proposition \ref{transfer}
motivates a crucial definition:
\begin{definition}
A {\bf $\ast$-finite} (pronounced ``star-finite'') subset of
$\mathbb{R}^\ast$ is one that can be placed in {\em internal}
one-to-one correspondence with $\set{1,\dots,n}$ for some
$n\in\mathbb{N}^\ast$.

Armed with this definition, readers not previously familiar
with nonstandard analysis can go back and read the description
of the example.
\end{definition}
\end{example}
\begin{example}\label{middlescale}
%
% In the next section this example and the two that follow it
% will afford counterexamples to some perhaps otherwise plausible
% propositions.
% 
This scale is not ``divided,'' but is ``divisible''
-- those terms will be defined in \S\ref{divisibility}.
Let $\mathbb{A}$ be the Boolean algebra of all subsets
of a set $\set{a,b,c}$, i.e., meet, join, complement are
the union, intersection, and complement operations on sets.
The set $\mathcal{R}=\set{0,\alpha,\beta,\gamma,\delta,1}$,
partially ordered as in Figure~\ref{midscalefig}, is a
scale.
   \input epsf
     \begin{figure}[htb]
     \begin{center}
     \leavevmode
     \hbox{%
     \epsfxsize=2.5in
     \epsfysize=2.5in
     \epsffile{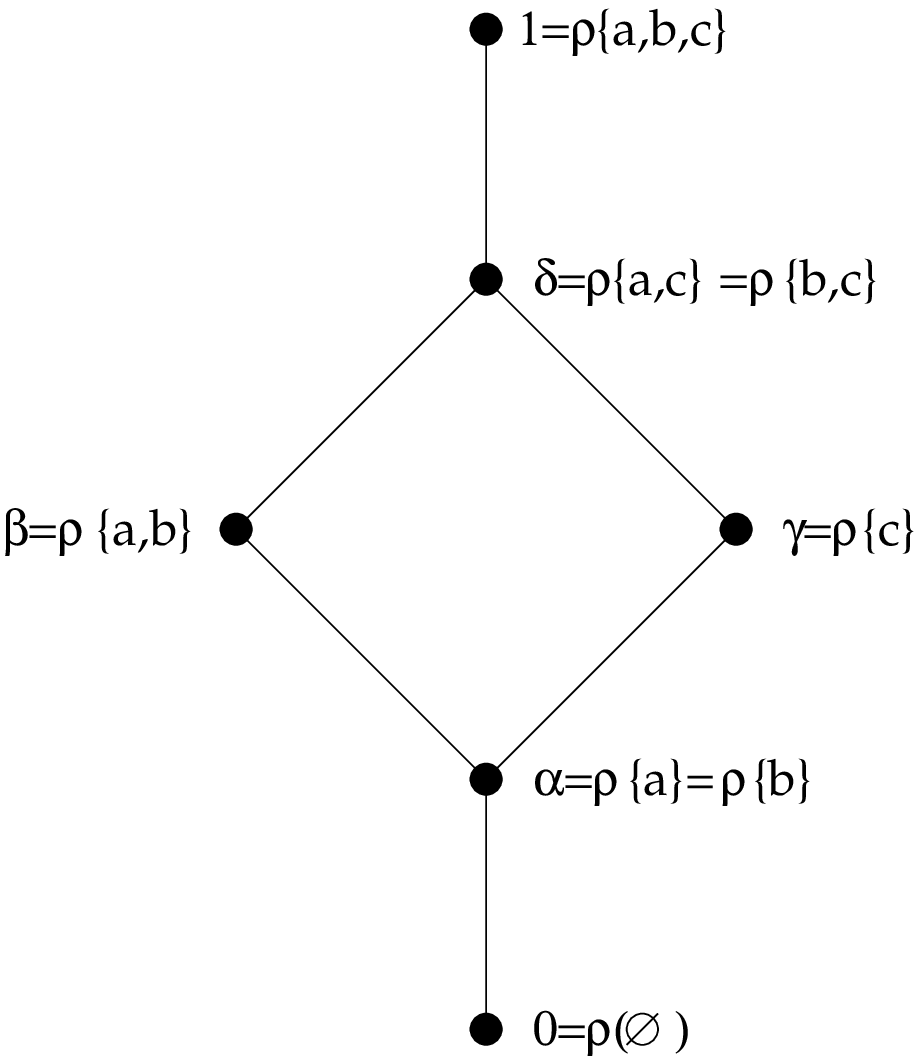}}
     \caption{}
     \label{midscalefig}
   \end{center}
 \end{figure}
\end{example}
\begin{example}
Extend the partial ordering of the previous example so that
$\beta<\gamma$, making $\mathcal{R}$ linearly ordered.
The same mapping into the same set, but with a different ordering
of that set, is a different scaling on the same Boolean algebra.
This scaling is isomorphic as a scaling to any measure $\mu$ on the
set of all subsets of $\set{a,b,c}$ satisfying
$$\mu\set{a}=\mu\set{b}
\quad\text{ and }\quad
\mu\set{a,b}<\mu\set{c}.$$
So measures that are not scalar multiples of each other
can be isomorphic to each other as scalings.
\end{example}
\begin{example}\label{quasicomplement}
By now Figure~\ref{quasicomplementfig} should be self-explanatory.
In \S\ref{relcomp}
this scale will provide us with an example of something that
``ought to be'' a relative complement but is not.
   \input epsf
     \begin{figure}[htb]
     \begin{center}
     \leavevmode
     \hbox{%
     \epsfxsize=3.0in
     \epsfysize=3.0in
     \epsffile{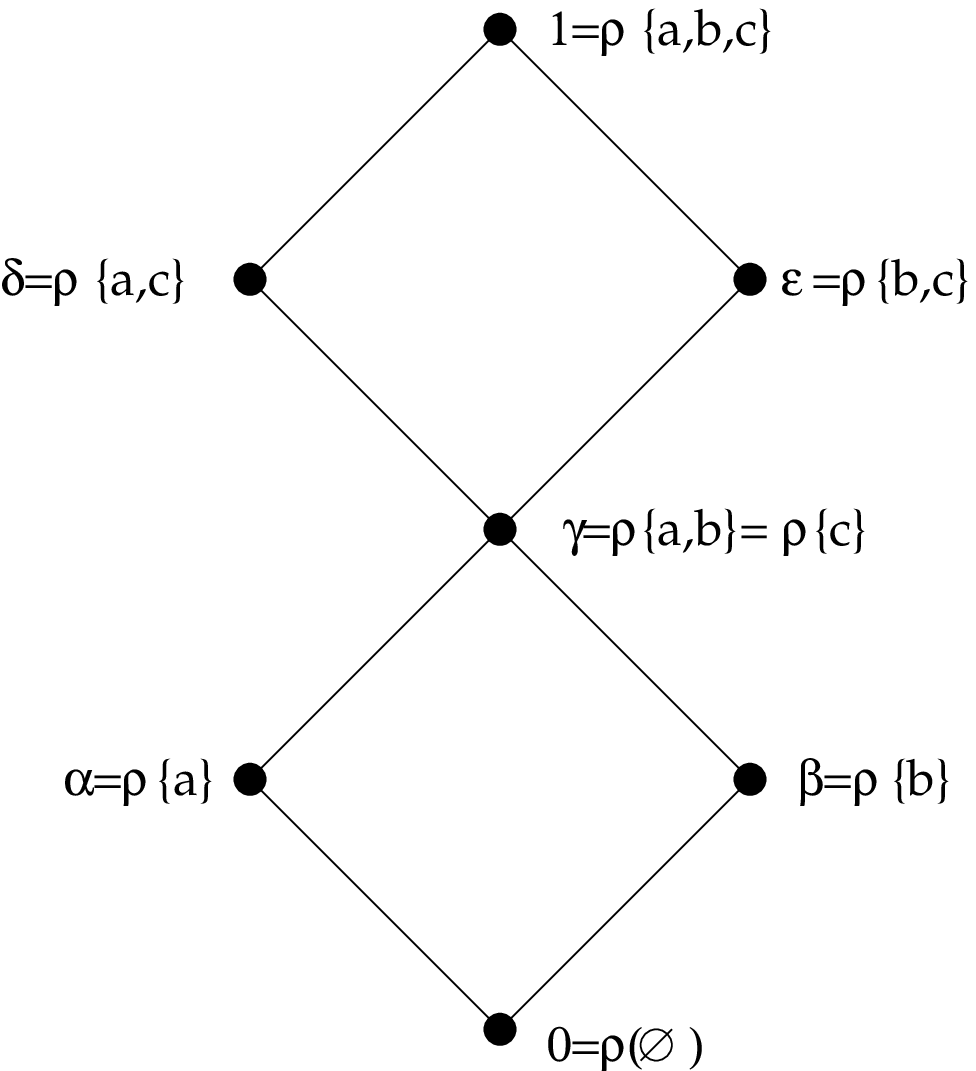}}
     \caption{}
     \label{quasicomplementfig}
   \end{center}
 \end{figure}
\end{example}
\begin{example}\label{balanced}
Let $\mathbb{A}$ be the Boolean algebra of all subsets
of $\set{a,b,c}$, and let $\mathcal{R}$ be the set of
all such subsets, partially ordered by saying that
$A<B$ iff $\left|A\right|<\left|B\right|$ (but note that
$\left|A\right|=\left|B\right|$ does not imply $\rho(A)=\rho(B)$).
See Figure~\ref{two-pointed}.
 \input epsf
  \begin{figure}[htb]
    \begin{center}
     \leavevmode
     \hbox{%
     \epsfxsize=2in
     \epsfysize=2.5in
     \epsffile{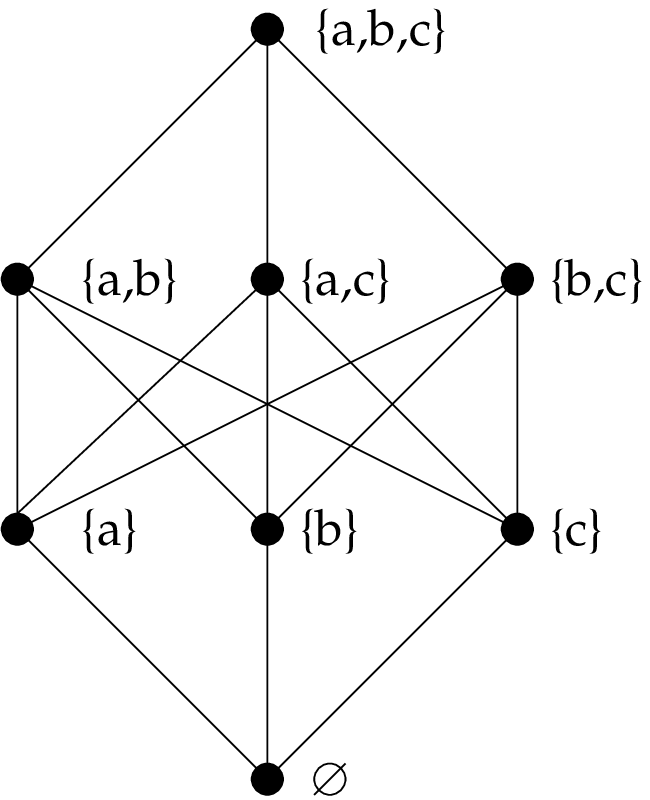}}
     \caption{}
     \label{two-pointed}
    \end{center}
  \end{figure}
Let $\rho$ map each subset of $\set{a,b,c}$ to itself.
Then $\rho$ is a scaling from $\mathbb{A}$ into $\mathcal{R}$.
One point of this example is that a scaling $\rho$ for which
$\rho(a)$, $\rho(b)$, and $\rho(c)$, are pairwise incomparable
can nonetheless have a properly more extensive partial ordering
than does the domain of the scaling.
Another point of this example is that this is another case in
which a scale's relative complementation is obviously not
determined by its partial ordering.  Finally, this is a finite
scale that is not a lattice; for example, there is no smallest
element that is $\geq$ both $\set{a}$ and $\set{b}$.
\end{example}

On the Boolean algebra of subsets of $\set{a,b,c}$, up to
isomorphism, there exist 17 one-to-one scalings, including
two whose images are linearly ordered scales, and 10 scalings
that are not one-to-one, including six whose images are linearly
ordered scales.  The last three examples above are of course
among these 27 scalings.  On the Boolean algebra of subsets of
$\set{a,b,c,d}$ there are 14 one-to-one scalings whose images
are linearly ordered, and many scalings that are less well-behaved.

\section{Additivity and its consequences}\label{addition}

\subsection{Addition}

The following lemma is an easy consequence of
Definition~\ref{BasicDefinitions}, but to get from this lemma
to Theorem~\ref{additivity}, the result that explains the title
of this section, is less straightforward.
\begin{lemma}\label{basiclemma}
If $\rho(x)\leq\rho(y)$ and $y \wedge z=0$,
then $\rho(x \vee z)\leq\rho(y \vee z)$,
and similarly if ``$<$'' replaces ``$\leq$'' throughout.
\end{lemma}
\begin{proof}
Let $x_1=x \wedge \ssim z \leq x$.
Then $x_1,y\in[0,\ssim z]$ and $\rho(x_1) \leq \rho(y)$.  Consequently
$\rho(\ssim x_{1 \ [0,\sim z]}) \geq \rho(\ssim y_{[0,\sim z]})$.
This reduces to
$\rho(\ssim z \wedge \ssim x_1) \geq \rho(\ssim z \wedge \ssim y)$,
whence we get $\rho(x_1 \vee z) \leq \rho(y \vee z)$.
Since $z \vee x_1=z \vee x$ we get
$\rho(x \vee z) \leq \rho(y \vee z)$.
For strict inequalities the proof is similar.
\end{proof}

The proof of the next result uses Lemma~\ref{basiclemma}
three times, but the three parts of the proof are not really
parallel to each other.
\begin{theorem}[Basic scalings are finitely additive]\label{additivity}
If $x \wedge y=0$ then $\rho(x \vee y)$ depends on $x$ and $y$
only through $\rho(x)$ and $\rho(y)$, and in a strictly increasing
fashion.  In other words, if $u \wedge v=0=x \wedge y$,
$\rho(u)=\rho(x)$, and $\rho(v)=\rho(y)$, then
$\rho(u \vee v)=\rho(x \vee y)$, and if ``$<$'' replaces
``$=$'' in either or both of the assumed equalities between
values of $\rho$, then ``$<$'' replaces ``$=$'' in the conclusion.
\end{theorem}
\begin{proof}
To prove ``$=$'' it suffices to prove both
``$\leq$'' and ``$\geq$''.  By symmetry we need only
do the first.
Although the proof that $\rho(u \vee v) \leq \rho(x \vee y)$
must rely on the fact that $\rho(u) \leq \rho(x)$ and
$\rho(v) \leq \rho(y)$, and that $x \wedge y=0$,
the assumption that $u \wedge v=0$ is needed
only for proving the inverse inequality ``$\geq$''.

Let $u_1=u \wedge \ssim y$, \ \ \ $y_1=y \wedge \ssim u$, \ \ \
$w=u \wedge y$.  We use Lemma~\ref{basiclemma} three times -- once
with $u_1$ in the role of $z$, once with $y_1$ in that role, and
once with $w$ in that role.
\begin{tabular}{c|c}
\begin{minipage}[t]{55mm}
By definition of $u_1$, we have $u_1 \wedge y=0$.
Therefore by Lemma~\ref{basiclemma} we can add $u_1$
to both sides of the inequality
\begin{center}
$\rho(v) \leq \rho(y)$
\end{center}
to get \ $\rho(u_1 \vee v) \leq \rho(u_1 \vee y)$. \vspace{2mm}
\end{minipage} &
\begin{minipage}[t]{53mm}
By hypothesis $x \wedge y=0$.  By definition of $y_1$, this
implies $x \wedge y_1=0$.  Therefore by Lemma~\ref{basiclemma}
we can add $y_1$ to both sides of the inequality
\begin{center}
$\rho(u) \leq \rho(x)$
\end{center}
to get \ $\rho(u \vee y_1) \leq \rho(x \vee y_1)$. \vspace{2mm}
\end{minipage}
\end{tabular}
Since $u_1 \vee y=u \vee y = u \vee y_1$, we get
\begin{equation}\label{third}
\rho(u_1 \vee v) \leq \rho(x \vee y_1).
\end{equation}
Before applying Lemma~\ref{basiclemma} the third time,
we must check that $(x \vee y_1) \wedge w=0$.
For this we need both the assumption that $x \wedge y=0$ and
the definitions of $y_1$ and $w$.
Then Lemma~\ref{basiclemma} applied to (\ref{third}) gives us
$$\rho(u_1 \vee v \vee w) \leq \rho(x \vee y_1 \vee w).$$
The definitions of $u_1$, $y_1$, and $w$, imply that
$u_1 \vee v \vee w=u \vee v$ and $x \vee y_1 \vee w=x \vee y$,
and the desired inequality follows.  Finally, use the statement
about strict inequalities in Lemma~\ref{basiclemma} to justify
the statement about strict inequalities in the conclusion of
the Theorem.
\end{proof}
\begin{corollary}\label{AdditionDefinition}
If $x \wedge y=0$, then $\rho(x)+\rho(y)$ can be defined
unambiguously as \mbox{$\rho(x \vee y)$}.
\end{corollary}

Addition is not everywhere-defined:
\begin{proposition}\label{partialoperation}
For $\zeta,\eta\in\mathcal{R}$, the sum $\zeta+\eta$
exists only if $\zeta\leq\ssim\eta$,
or, equivalently, $\eta\leq\ssim\zeta$.
\end{proposition}
\begin{proof}
If $\zeta+\eta$ exists then for some $x,y\in\mathbb{A}$
we have $\rho(x)=\zeta$, $\rho(y)=\eta$, and $x\wedge y=0$.
But $x\wedge y=0$ is equivalent to $x\leq\ssim y$,
and that implies $\zeta=\rho(x)\leq\rho(\ssim y)=\ssim\eta$.
\end{proof}

A sum $\zeta+\eta$ may be undefined even when
$\zeta\leq\ssim \eta$, simply because there are
no two elements $x,y \in \mathbb{A}$ such that
$x \wedge y=0$, $\rho(x)=\zeta$, and $\rho(y)=\eta$.
The remedy to this unpleasant situation is in \S\ref{divisibility}.
The problem occurs in Example~\ref{middlescale}, where
$\alpha\leq\ssim\beta$, but $\alpha+\beta$ is nonetheless
undefined.  The addition table for that example appears
in Figure~\ref{additiontable}.
\begin{figure}[htb]
   \begin{center}
$$\begin{array}{|c|cccccc|}
\hline
\ + \    & 0 \quad & \alpha \quad & \beta \quad &
                     \gamma \quad & \delta \quad & 1 \quad \\
\hline
\ 0 \ & 0 \quad & \alpha \quad & \beta \quad & \gamma \quad &
                          \delta \quad & 1 \quad \\
\ \alpha \ & \alpha \quad & \beta \quad & \text{\fbox{?}} \quad &
               \delta \quad & 1 \quad & \boxtimes \quad \\
\ \beta \ & \beta \quad & \text{\fbox{?}} \quad &
      \boxtimes \quad & 1 \quad & \boxtimes \quad &
               \boxtimes \quad \\
\ \gamma \ & \gamma \quad & \delta \quad & 1 \quad &
     \boxtimes \quad & \boxtimes \quad & \boxtimes \quad \\ 
\ \delta \ & \delta \quad & 1 \quad &
    \boxtimes \quad & \boxtimes \quad & \boxtimes \quad &
                           \boxtimes \quad \\
\ 1 \ & 1 \quad & \boxtimes \quad & \boxtimes \quad & \boxtimes \quad &
     \boxtimes \quad & \boxtimes \quad \\
\hline
\end{array}$$
\caption{Addition table for Example~\ref{middlescale}.}
\label{additiontable}\label{1stuse}
   \end{center}
\end{figure}
A ``$\boxtimes$'' marks the places where
Proposition~\ref{partialoperation} explains why the entry
is undefined.  A ``\fbox{?}'' marks the other places where
the entry is undefined.  Example~\ref{nonarch} is a scale
for which this particular pathology -- that $\zeta+\eta$
may be undefined even though $\zeta\leq\ssim\eta$ -- never
occurs.

Some ways in which addition is obviously well-behaved are
these: For $\zeta,\eta\in\mathcal{R}$ we have
$\zeta+0=\zeta$, $\zeta+\ssim\zeta=1$, and
$\zeta+\eta=\eta+\zeta$, the existence of either of these
sums entailing that of the other.

What about associativity?  If $\zeta+(\eta+\theta)$ exists
then some $y,z\in\mathbb{A}$ whose images under $\rho$ are
$\eta$ and $\theta$ respectively, satisfy $y\wedge z=0$,
and some $x,w\in\mathbb{A}$ whose images under $\rho$ are
$\zeta$ and $\eta+\theta=\rho(y\vee z)$ satisfy $x\wedge w=0$.
Neither $y$ nor $z$ was assumed disjoint from $x$.
Can $w$ be split into disjoint parts whose images under
$\rho$ are those of $y$ and $z$?  Not always.  When
$\zeta+(\eta+\theta)$ and $(\zeta+\eta)+\theta$ both exist
are they always equal?  A partial answer is obvious:
\begin{proposition}\label{associativity}
If $x_1,\dots,x_n\in\mathbb{A}$
are {\em pairwise} disjoint and $\rho(x_i)=\zeta_i$ for
$i=1,\dots,n$, then $\zeta_1+\cdots+\zeta_n$ is unambiguously
defined.
\end{proposition}

\subsection{Duality, modularity, subtraction, relative
complementation, and de Morganism}

\subsubsection{Duality}
The array~(\ref{booleIDs}) of identities on
page~\pageref{booleIDs} defining the concept of Boolean algebra
has an evident symmetry: Interchange the roles of ``$\wedge$''
and ``$\vee$'' and of ``0'' and ``1'', and the identities in
that table are merely permuted among themselves.
If the partial ordering $\leq$ on a Boolean algebra
is regarded as part of the structure, interchange it with its
inverse $\geq$.  All consequences of those identities then
remain true if this same interchange of relations and operations
is applied to them.  The interchange leaves the operation of
complementation unchanged, i.e., that operation is self-dual.
That much is well-known.  The same thing applies not only to
Boolean algebras but also to scales generally.
In particular, the dual
%
%   Proposition~\ref{partialoperation} implies that if
%   $\rho(x)\not\geq\ssim\rho(y)$, or equivalently
%   $\rho(y)\not\geq\ssim\rho(x)$, then $\rho(x \vee y)<1$.
%   This is a weak sort of finite sub-additivity whose statement
%   does not rely on the concept of addition and whose proof
%   requires only absolute complements and not relative complements
%   (relative complements are used in the proof of Lemma~\ref{basiclemma}).
%
%   More important for our purposes is the dual of
%
of Theorem~\ref{additivity} says:
\begin{quote}
{\em If $x\vee y=1$ then $\rho(x\wedge y)$ depends on $x$
and $y$ only through $\rho(x)$ and $\rho(y)$, and in a
strictly increasing fashion.}
\end{quote}
Therefore a ``dual-addition'' is unambiguously defined.
We shall call the values of this operation ``dual-sums''
and write
$$\rho(x \wedge y)=\rho(x)\oplus\rho(y) \text{ if } x \vee y=1.$$
The dual of Proposition~\ref{partialoperation} says that
$\zeta\oplus\eta$ exists only if $\zeta\geq\ssim\eta$, or
equivalently, $\eta\geq\ssim\zeta$.
The dual-addition table for Example~\ref{middlescale} on
page~\pageref{middlescale} is constructed by first reflecting
the interior, but not the margins, of the table in
Figure~\ref{additiontable} about the diagonal that contains
only 1's, and then replacing each entry in the interior, but
not in the margins, by its complement:\label{2nduse}
$0 \leftrightarrow 1$, $\alpha \leftrightarrow \delta$,
$\beta \leftrightarrow \gamma$, $\boxtimes\leftrightarrow\boxtimes$,
$\text{\fbox{?}}\leftrightarrow\text{\fbox{?}}$.

\subsubsection{Modularity}
The following ``modular law'' is the essential
tool for dealing with subtraction, cancellation, and
relative complementation.
\begin{lemma}\label{modularlaw}
If $x\wedge y=0$ and $y\vee z=1$
then no ambiguity comes from writing
$$\rho(x)+\rho(y)\oplus\rho(z).$$
In other words, we have
$\left\{\rho(x)+\rho(y)\right\}\oplus\rho(z)
=\rho(x)+\left\{\rho(y)\oplus\rho(z)\right\}$.
In particular, the sums and dual-sums exist.
\end{lemma}
\begin{proof}
Since $x \wedge y=0$ we also have $x \wedge (y \wedge z)=0$
so $\rho(x)+\left\{\rho(y)\oplus\rho(z)\right\}$
exists and is equal to $\rho(x \vee (y \wedge z))$.
Since $y \vee z=1$ we have $(x \vee y)\vee z=1$, so
$\left\{\rho(x)+\rho(y)\right\}\oplus\rho(z)$
exists and is equal to $\rho((x \vee y)\wedge z)$.
Finally, the two identities $x \wedge y=0$ and $y \vee z=1$
entail that $x \vee (y \wedge z)=(x \vee y)\wedge z$.
\end{proof}
\newpage
\subsubsection{Subtraction}
\begin{proposition}\label{subtraction}
The functions
$$
\zeta\mapsto\zeta+\eta \quad\text{ and }\quad
\zeta\mapsto\zeta\oplus\ssim\eta
$$
are inverses.  In particular, the not-everywhere-defined
nature of the operations involved does not prevent the image
of each of these functions from coinciding with the domain
of the other.
\end{proposition}
In other words, subtraction of $\eta$ from $\zeta$
yields $\zeta-\eta=\zeta\oplus\ssim\eta$.
\begin{proof}
Suppose $x\wedge y=0$, $\rho(x)=\zeta$, and $\rho(y)=\eta$.
Since $y\vee\ssim y=1$, Lemma~\ref{modularlaw} (the modular law)
applies:
$$
(\zeta+\eta)\oplus\ssim\eta=\zeta+(\eta\oplus\ssim\eta)
=\zeta+0=\zeta.
$$
So the second function is a left-inverse of the first.
To prove the first is a left-inverse of the second,
dualize, interchanging ``$\wedge$'' with ``$\vee$'',
``0'' with ``1'', and ``$+$'' with ``$\oplus$''.
\end{proof}
\begin{proposition}
The difference $\zeta-\eta$ exists only if $\zeta\geq\eta$.
\end{proposition}
\begin{proof}
The dual of Proposition~\ref{partialoperation} implies
$\zeta\oplus\ssim\eta$ exists only if $\zeta\geq\eta$.
\end{proof}
But $\zeta-\eta$ is sometimes undefined even when $\zeta\geq\eta$.
In Example~\ref{middlescale}, we have $\gamma>\alpha$,
but no members $x,y$ of the domain simultaneously satisfy
$y>x$, $\rho(y)=\gamma$, and $\rho(x)=\alpha$.  Thus
we cannot subtract $\alpha$ from $\gamma$.  This difficulty
will be remedied in \S\ref{divisibility}.
\subsubsection{Relative complementation}\label{relcomp}
Lemma~\ref{modularlaw} (the modular law) can be used to prove
that a scaling induces an operation of additive relative
complementation on its image.  Like addition and subtraction,
this is not every\-where-defined.
\begin{proposition}
If $x\in [a,b]$ then $\rho(\ssim x_{[a,b]})$ depends on $x$,
$a$, and $b$ only through $\rho(x)$, $\rho(a)$, and $\rho(b)$.
\end{proposition}
Recall the difficulty: In effect the proposition says
if $x\in[a,b]$, $y \in [c,d]$, and
$\rho(a)=\rho(c) \leq \rho(x)=\rho(y) \leq \rho(b)=\rho(d)$
then $\rho(\ssim x_{[a,b]})=\rho(\ssim y_{[c,d]})$.
But only in case $[a,b]=[c,d]$ is this immediate
from Definition~\ref{BasicDefinitions}.
\begin{proof}
\enlargethispage{1cm}
Since $a\leq x\leq b$ we have $a\wedge\ssim x=0$ and
$(\ssim x)\vee b=1$.  Therefore, by Lemma~\ref{modularlaw}
(the modular law), the following is unambiguously defined.
$$
\rho(\ssim x_{[a,b]})=\rho(a)+\rho(\ssim x)\oplus\rho(b).
$$
\end{proof}
Example~\ref{middlescale} shows why additive relative
complementation on a scale is not\label{3rduse}
every\-where-defined.  In that example $\alpha$ is its
own additive complement relative to the interval $[0,\beta]$.
But the additive complement of $\alpha$ relative
to the interval $[0,\gamma]$ does not exist even though
$0<\alpha<\gamma$, because there do not exist
$0<x<y$ in the domain of $\rho$ whose respective
images under $\rho$ are $0$, $\alpha$, and $\gamma$.
This will be remedied in \S\ref{divisibility}.

In \S\ref{examples} we remarked that Example~\ref{quasicomplement}
``provide[s] us with an example of something that
`ought to be' a relative complement but is not.''
In that example we have $\alpha<\gamma<\varepsilon$ and
$\beta<\gamma<\delta$, but no complements
$\ssim\gamma_{[\alpha,\varepsilon]}$ or
$\ssim\gamma_{[\beta,\delta]}$ exist, even though
$\alpha+(\ssim\gamma)\oplus\varepsilon=\gamma$
and $\beta+(\ssim\gamma)\oplus\delta=\gamma$ do exist.

The operation of additive relative complementation on a
scale depends not only on the partial ordering of the scale
but also on the scaling.
This can be seen by considering Example~\ref{balanced}.
The additive complement of $\set{a,b}$ relative to the interval
from $\set{a}$ to $\set{a,b,c}$ is $\set{a,c}$, not only
in the domain, but also in the range!  But nothing in the
partial ordering of that scale makes $\set{a,c}$ a better
candidate than $\set{b,c}$ to be the additive relative complement.
Rather, it is singled out by the ordering of the domain.
\subsubsection{de Morganism}
The next proposition is immediate
from the results of this section and will be useful
in \S\ref{divisibility}:
\begin{proposition}\label{demorgan}
\begin{eqnarray*}
\text{If }\zeta+\eta\text{ exists, then so does }
(\ssim\zeta)\oplus(\ssim\eta)\text{, and }
\ssim(\zeta+\eta) & = & (\ssim\zeta)\oplus(\ssim\eta). \\
\text{If }\zeta\oplus\eta\text{ exists, then so does }
(\ssim\zeta)+(\ssim\eta)\text{, and }
\ssim(\zeta\oplus\eta) & = & (\ssim\zeta)+(\ssim\eta).
\end{eqnarray*}
\end{proposition}
\subsubsection{Technical lemma on inequalities}
The following lemma will be useful in \S\ref{divisibility}.
\begin{lemma}\label{InequalityLemma}
\begin{eqnarray}
\text{If }\zeta\leq\eta\text{ and }\zeta+\theta\text{ and }
\eta+\theta\text{ exist, then }
\zeta+\theta\leq\eta+\theta, \label{InequalityLemma1} \\
\text{If }\zeta\leq\eta\text{ and }\zeta\oplus\theta\text{ and }
\eta\oplus\theta\text{ exist, then }
\zeta\oplus\theta\leq\eta\oplus\theta, \label{InequalityLemma2} \\
\text{If }\zeta\leq\eta\text{ and }\zeta-\theta\text{ and }
\eta-\theta\text{ exist, then }
\zeta-\theta\leq\eta-\theta, \label{InequalityLemma3}
\end{eqnarray}
and all three statements remain true if ``$<$'' replaces
both occurrences of ``$\leq$''.
\end{lemma}
\begin{proof}
(\ref{InequalityLemma1}) follows from the conjunction of
Theorem~\ref{additivity} with the definition embodied in
Corollary~\ref{AdditionDefinition}.
(\ref{InequalityLemma2}) is the dual of (\ref{InequalityLemma1}).
(\ref{InequalityLemma3}) follows from the conjunction of
(\ref{InequalityLemma2}) with Proposition~\ref{subtraction}
and the strictly decreasing nature of (absolute) complementation.
\end{proof}
\subsection{Measures on scales}\label{measures.too}
\begin{definition}
A {\bf finitely additive measure}
on a scale $\mathcal{R}$ is a strictly increasing mapping
$\mu:\mathcal{R}\rightarrow[0,\infty)\subseteq\mathbb{R}$
satisfying
\begin{enumerate}
\item
For $\zeta\in\mathcal{R}$, if $\zeta>0$ then $\mu(\zeta)>0$, and
\item
For $\zeta,\eta\in\mathcal{R}$, if $\zeta+\eta$ exists
then $\mu(\zeta+\eta)=\mu(\zeta)+\mu(\eta)$.
\end{enumerate}
A scale is {\bf measurable} if it is the domain of a measure.
\end{definition}
The words ``strictly increasing'' would be redundant if the
domain $\mathcal{R}$ were a Boolean algebra.  They are also
redundant in what we shall call ``divided'' scales -- to be
defined in the next section.  That they are not redundant in
this more general setting is shown by this example:
Let $\rho$ be defined on the set of all four subsets
of $\set{a,b}$ and suppose $0<\rho(\set{a})<\rho(\set{b})<1$.
Let $0=\mu(0)<\mu(\rho(\set{a})=2/3\not\leq 1/3=\mu(\rho(\set{b}))
<\mu(1)=1$.

Clearly this generalizes Definition~\ref{measures.on.BAs}.
Moreover, if $\mu:\mathcal{R}\rightarrow[0,\infty)$ is a
measure and $\rho:\mathbb{A}\rightarrow\mathcal{R}$ is the
basic scaling that induces the scale-structure on $\mathcal{R}$,
then $\mu\circ\rho:\mathbb{A}\rightarrow[0,\infty)$ is a
measure on the Boolean algebra $\mathbb{A}$.

\section{Divisibility and measurability}\label{divisibility}

\subsection{Dividedness}\label{dividedness}

In \S~\ref{addition} we saw three pathologies:
\begin{enumerate}
\item
Although $\zeta\leq\ssim\eta$, or equivalently $\eta\leq\ssim\zeta$,
is necessary for the existence of $\zeta+\eta$, in some scales it
is not sufficient because there may be no $x,y\in\mathbb{A}$ for
which $x\wedge y=0$ and $\rho(x)=\zeta$ and $\rho(y)=\eta$;
\item
Although $\zeta\geq\eta$ is necessary of the existence of
$\zeta-\eta$, in some scales it is not sufficient because there
may be no $x,y\in\mathbb{A}$ for which $x>y$ and $\rho(x)=\zeta$
and $\rho(y)=\eta$;
\item
Although $\zeta\leq\eta\leq\theta$ is necessary for the existence
of $\ssim\eta_{[\zeta,\theta]}$, in some scales it is not sufficient
because there may be no $x,y,z\in\mathbb{A}$ for which $x\leq y\leq z$
and $\rho(x)=\zeta$, $\rho(y)=\eta$, and $\rho(z)=\theta$.
\end{enumerate}
\begin{proposition}
These three pathologies are equivalent, i.e., in any scale
in which one of them occurs, so do the others.
\end{proposition}
\begin{proof}
The sum $\zeta+\eta$ and the difference
$(\ssim\eta)-\zeta=(\ssim\eta)\oplus(\ssim\zeta)$
are complements of each other, and complementation on a scale
is a bijection.  This suffices for equivalence of (1) and (2).
Existence of the relative complement $\ssim\eta_{[\zeta,\theta]}$
is equivalent to the existence of both the sum $\zeta+(\ssim\eta)$
and the difference $\theta-\eta$.  This suffices for equivalence of
(3) with its predecessors.
\end{proof}
A simple law additional to those that define a scale is the remedy.
\begin{definition}
A basic scaling $\rho:\mathbb{A}\rightarrow\mathcal{R}$ is
{\bf divided} if whenever $\rho(x)<\rho(y)$ then some
$y_1<y$ satisfies $\rho(y_1)=\rho(x)$.  The domain
$\mathbb{A}$ and the range $\mathcal{R}$ will also be
called ``divided'' if $\rho$ is divided.
\end{definition}
In case $\mathbb{A}$ is an algebra of subsets of a set and $\rho$
is a measure, this says the measurable set $y$ is the union of
smaller sets, one of which has the same measure as $x$.
The reader can check that Examples~\ref{isomorphism},
\ref{nonarch}, and~\ref{fullnonarch} are divided, and
Example~\ref{middlescale} is not divided.
\begin{theorem}\label{divided-defined}
If a scale $\mathcal{R}$ is divided, then
for $\zeta,\eta,\theta\in\mathcal{R}$,
  \begin{enumerate}
    \item\label{divided-defined1}
      $\phantom{\text{T}}
      \zeta+\eta$ exists if $\zeta\leq\ssim\eta$
      (or equivalently, if $\eta\leq\ssim\zeta$).
    \item\label{divided-defined2}
      $\phantom{\text{T}}
      \zeta-\eta$ exists if $\zeta\geq\eta$.
    \item\label{divided-defined3}
      $\phantom{\text{T}}
      \ssim\eta_{[\zeta,\theta]}$ exists if
      $\zeta\leq\eta\leq\theta$.
  \end{enumerate}
\end{theorem}
\begin{proof}
Item (2) follows immediately from the definition
of dividedness.

It follows that under the assumptions of item (1),
(i.e., that $\zeta\leq\ssim\eta$), $\zeta$ can
be subtracted from $\ssim\eta$.  Then we have
\mbox{$(\ssim\eta)-\zeta=(\ssim\eta)\oplus(\ssim\zeta)$};
in particular, this latter dual-sum exists.  By
Proposition~\ref{demorgan}, on ``de-Morganism,''
so does the sum $\zeta+\eta$.

Item (3) follows from the conjunction of items (1) and (2)
and the observation that
$\ssim\eta_{[\zeta,\theta]}=\zeta+(\theta-\eta)$.
(Note that the parentheses in ``$\zeta+(\theta-\eta)$''
need to be where they are.)
\end{proof}

\subsection{Atoms}

\begin{definition}
  \begin{enumerate}
    \item An element $x\neq 0$ of a Boolean algebra
$\mathbb{A}$ is an {\bf atom} of $\mathbb{A}$ if the
interval $[0,x]$ contains only 0 and $x$.
    \item A Boolean algebra $\mathbb{A}$ is {\bf atomic} if
for every $y\in\mathbb{A}$ there is some atom
$x\in\mathbb{A}$ such that $x\leq y$.
    \item A Boolean algebra $\mathbb{A}$ is {\bf atomless}
if it contains no atoms.
\end{enumerate}
\end{definition}

\begin{example}
The Boolean algebra of all subsets of a set is atomic.  Each
singleton, i.e., each subset with only one member, is an atom.
\end{example}

\begin{example}
The Boolean algebra of all clopen (i.e., simultaneously
closed and open) subsets of the Cantor set is atomless.
\end{example}

\begin{example}
Adjoin a finite set of isolated points to the Cantor set.
The Boolean algebra of all clopen subsets of the resulting
space is neither atomic nor atomless.  The singleton of
each isolated point is an atom.
\end{example}

\subsection{Divisions}

\S\ref{dividedness} may appear to be suggesting that in
cases in which $\zeta+\eta$ does not exist, we should seek
some larger Boolean algebra and a correspondingly larger
scale in which the sum $\zeta+\eta$ will be found.  Hence
we have the following definition.
\begin{definition}
A {\bf division} of a basic scaling $\rho:\mathbb{A}\rightarrow\mathcal{R}$
is a scaling $\rho_1:\mathbb{A}_1\rightarrow\mathcal{R}_1$ such that
   \begin{enumerate}
     \item \ \ \ $\rho_1$ is divided;
     \item \ \ \ $\mathbb{A}$ is a subalgebra of $\mathbb{A}_1$;
     \item \ \ \ $\mathcal{R}$ is a sub-poset of $\mathcal{R}_1$;
     \item \ \ \ $\rho$ is the restriction of $\rho_1$ to $\mathcal{R}$; and
     \item \ \ \ No pair intermediate between
       $\left(\mathbb{A},\mathcal{R}\right)$ and
       $\left(\mathbb{A}_1,\mathcal{R}_1\right)$ satisfies 1-4.
     \end{enumerate}
The domain $\mathbb{A}_1$ and the range $\mathcal{R}_1$
of $\rho_1$ will also be called ``divisions'' of $\mathbb{A}$
and $\mathcal{R}$ respectively.
\end{definition}
\begin{example}
Regard the Boolean algebra $\mathbb{A}$ of all four subsets
of $\set{a,b}$ as a subalgebra of the Boolean algebra
$\mathbb{A}_1$ of all subsets of $\set{a,b_1,b_2}$ by
identifying $b$ with $\set{b_1,b_2}$, so that the atom $b$
has been split.  Then the scale on the right in
Figure~\ref{divisionfig1} is a division of the one on the left.
   \input epsf
     \begin{figure}[htb]
     \begin{center}
     \leavevmode
     \hbox{%
     \epsfxsize=3.6in
     \epsfysize=2.0in
     \epsffile{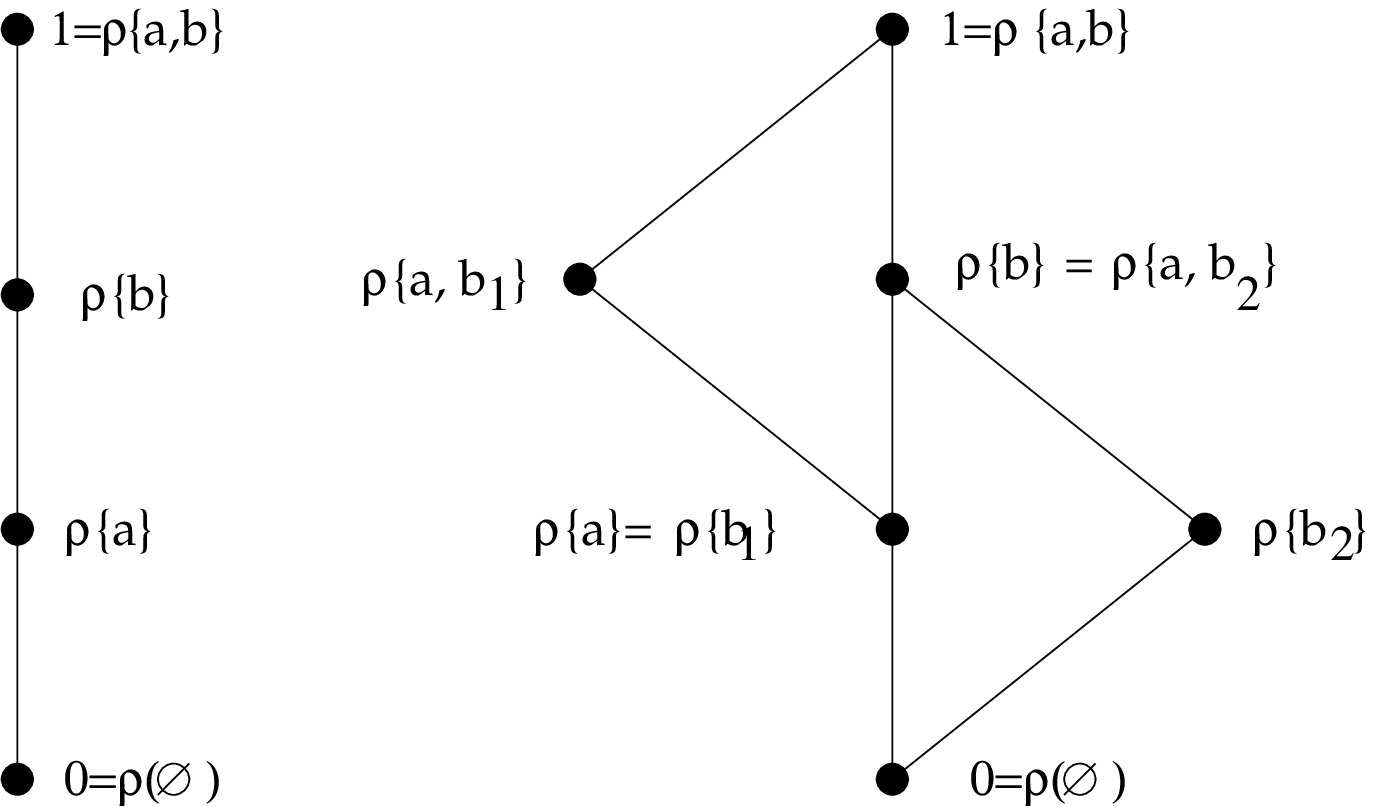}}
     \caption{}
     \label{divisionfig1}
     \end{center}
     \end{figure}
\end{example}
\begin{example}
Regard the Boolean algebra $\mathbb{A}$ of all eight subsets
of $\set{a,b,c}$ as a subalgebra of the Boolean algebra
$\mathbb{A}_1$ of all subsets of $\set{a,b,c_1,c_2}$ by
identifying $c$ with $\set{c_1,c_2}$, so that the atom $c$
has been split.
   \input epsf
     \begin{figure}[htb]
     \begin{center}
     \leavevmode
     \hbox{%
     \epsfxsize=2.5in
     \epsfysize=2.5in
     \epsffile{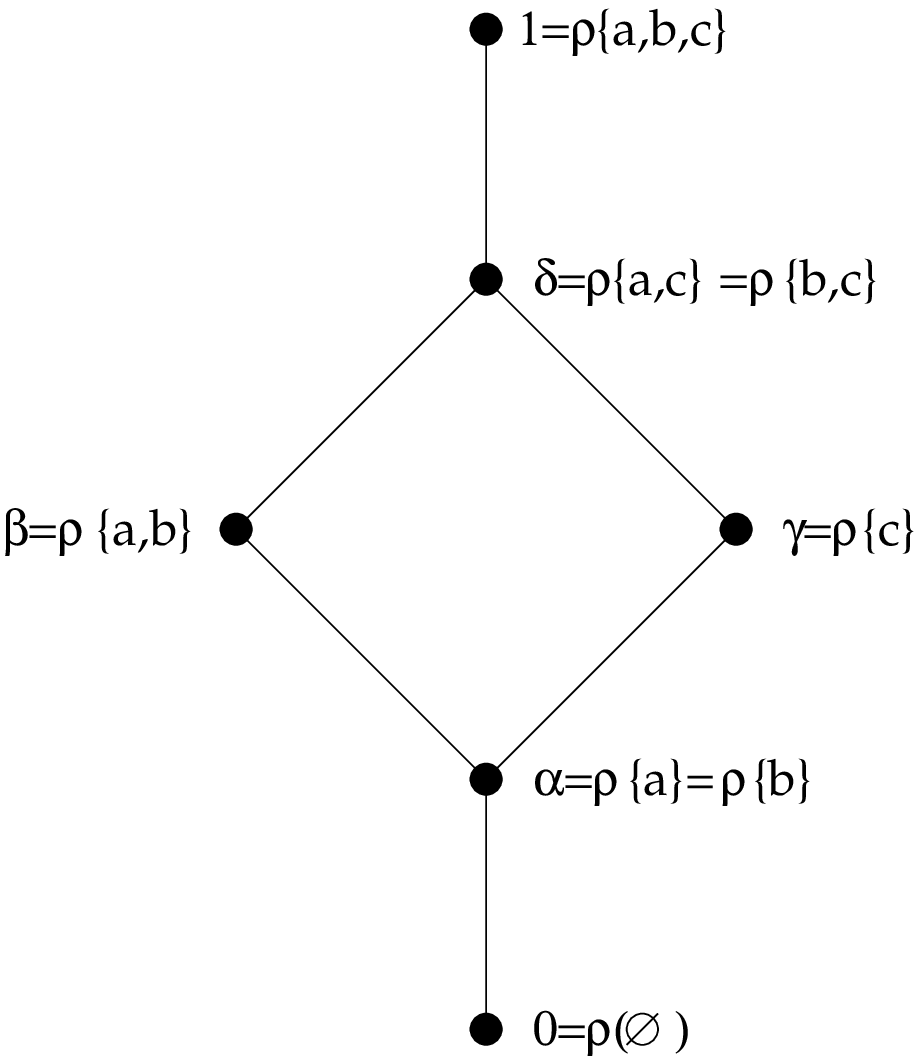}}
     \caption{}
     \label{BB}
     \end{center}
     \end{figure}
   \input epsf
     \begin{figure}[htb]
     \begin{center}
     \leavevmode
     \hbox{%
     \epsfxsize=4.5in
     \epsfysize=2.8in
     \epsffile{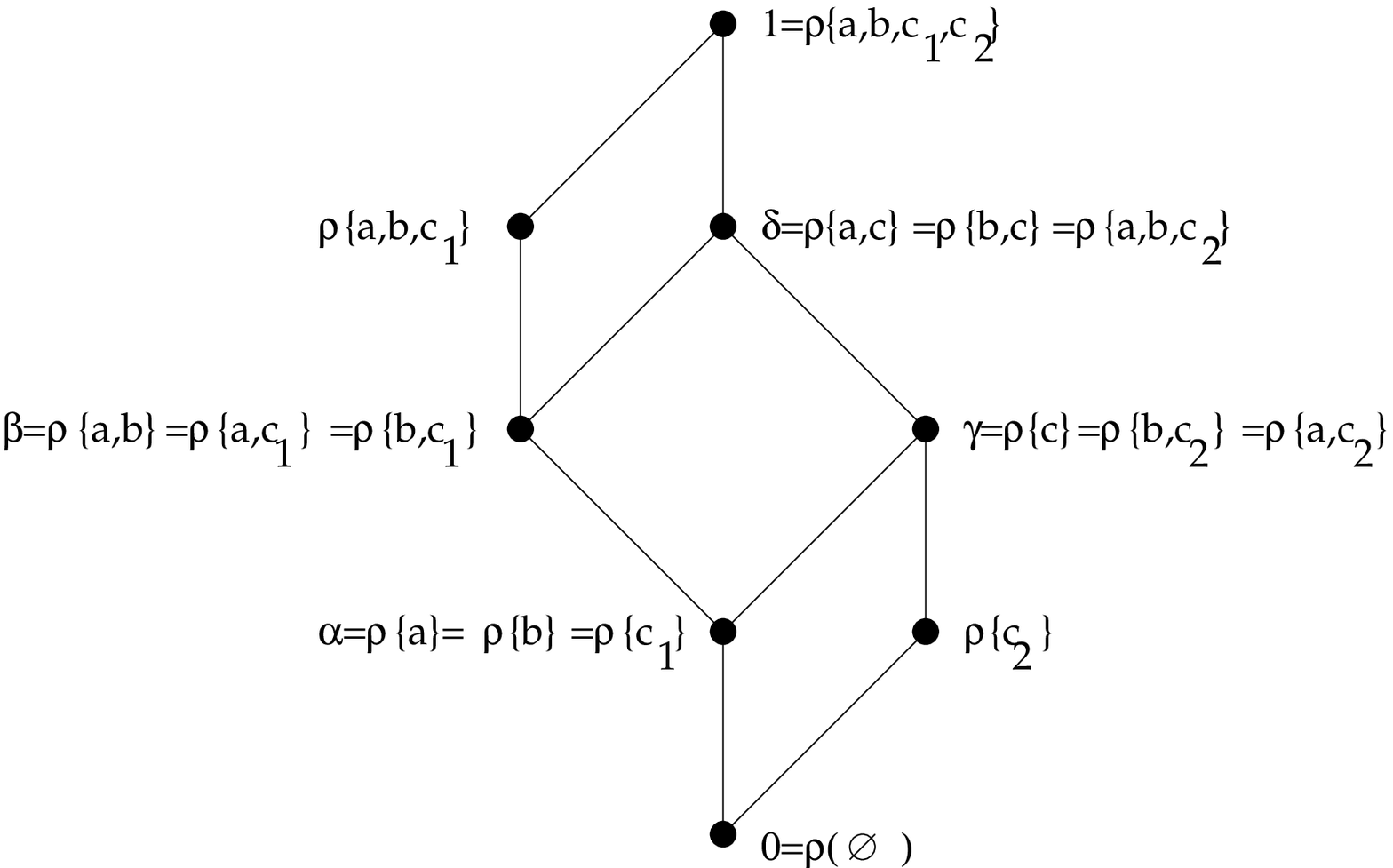}}
     \caption{}
     \label{AA}
     \end{center}
     \end{figure}
A division of the scale in Figure~\ref{BB}
appears in Figure~\ref{AA}.
\end{example}
\begin{example}\label{SplittingAlgorithm}
Suppose $\mathbb{A}$ is the Boolean algebra of all subsets
of $\set{a,b,c,d,e}$, and $\rho$ is the scaling arising
in the manner described in Example~\ref{convexscale} from
the convex set $\mathcal{C}$ of all measures $m$ on $\mathbb{A}$
that satisfy $m(\set{a,b,c})<m(\set{d,e})$.  Then
$\rho(\set{a,b,c})<\rho(\set{d,e})$, and the only other
sets $S,T\in\mathbb{A}$ for which $\rho(S)\leq\rho(T)$
are those for which $S\subseteq T$.  Split $a$ into
disjoint parts $a_d$ and $a_e$, similarly $b$ into
$b_d$ and $b_e$, and $c$ into $c_d$ and $c_e$.
Split $d$ into four disjoint parts $d_a$, $d_b$, $d_c$,
$d_{\text{ceterus}}$, and $e$ into $e_a$, $e_b$, $e_c$,
$e_{\text{ceterus}}$.  The convex set $\mathcal{C}$ is
then naturally identified with the set of all measures $m_1$
on the Boolean algebra $\mathbb{A}_1$ of all subsets of
$$
\set{a_d,a_e,b_d,b_e,c_d,c_e,d_a,d_b,d_c,d_{\text{ceterus}},
e_a,e_b,e_c,e_{\text{ceterus}}}
$$
that satisfy
$$
m_1(\set{a_d,a_e,b_d,b_e,c_d,c_e})
<m_1(\set{d_a,d_b,d_c,d_{\text{ceterus}},
e_a,e_b,e_c,e_{\text{ceterus}}}).
$$
Let $\mathcal{C}_1\subseteq\mathcal{C}$ be the smaller class
of measures $m_1$ that satisfy this inequality and also
$m_1(\set{a_d})=m_1(\set{d_a})$, $m_1(\set{a_e})=m_1(\set{e_a})$,
and so on.  (If we had had $\rho(\set{a,b,c})=\rho(\set{d,e})$,
i.e., ``$=$'' instead of ``$<$'', then we would have omitted
$d_{\text{ceterus}}$ and $e_{\text{ceterus}}$, which are
slack components.)  Then the scaling $\rho_1$ arising from
$\mathcal{C}_1$ in the manner of Example~\ref{convexscale}
is a division of $\rho$.
\end{example}

\subsection{Divisibility}

\subsubsection{Defined}

Does every finite scale have a division?
The next definition foreshadows the answer.
\begin{definition}
A scale is {\bf divisible} if it has a division;
otherwise it is {\bf indivisible}.
\end{definition}

\subsubsection{Adaptation of the Kraft-Pratt-Seidenberg counterexample}

\begin{example}\label{IndivisibleScale}
Figure~\ref{kraftposet} depicts
   \input epsf
     \begin{figure}[htb]
     \begin{center}
     \leavevmode
     \hbox{%
     \epsfxsize=4.5in
     \epsfysize=5.0in
     \epsffile{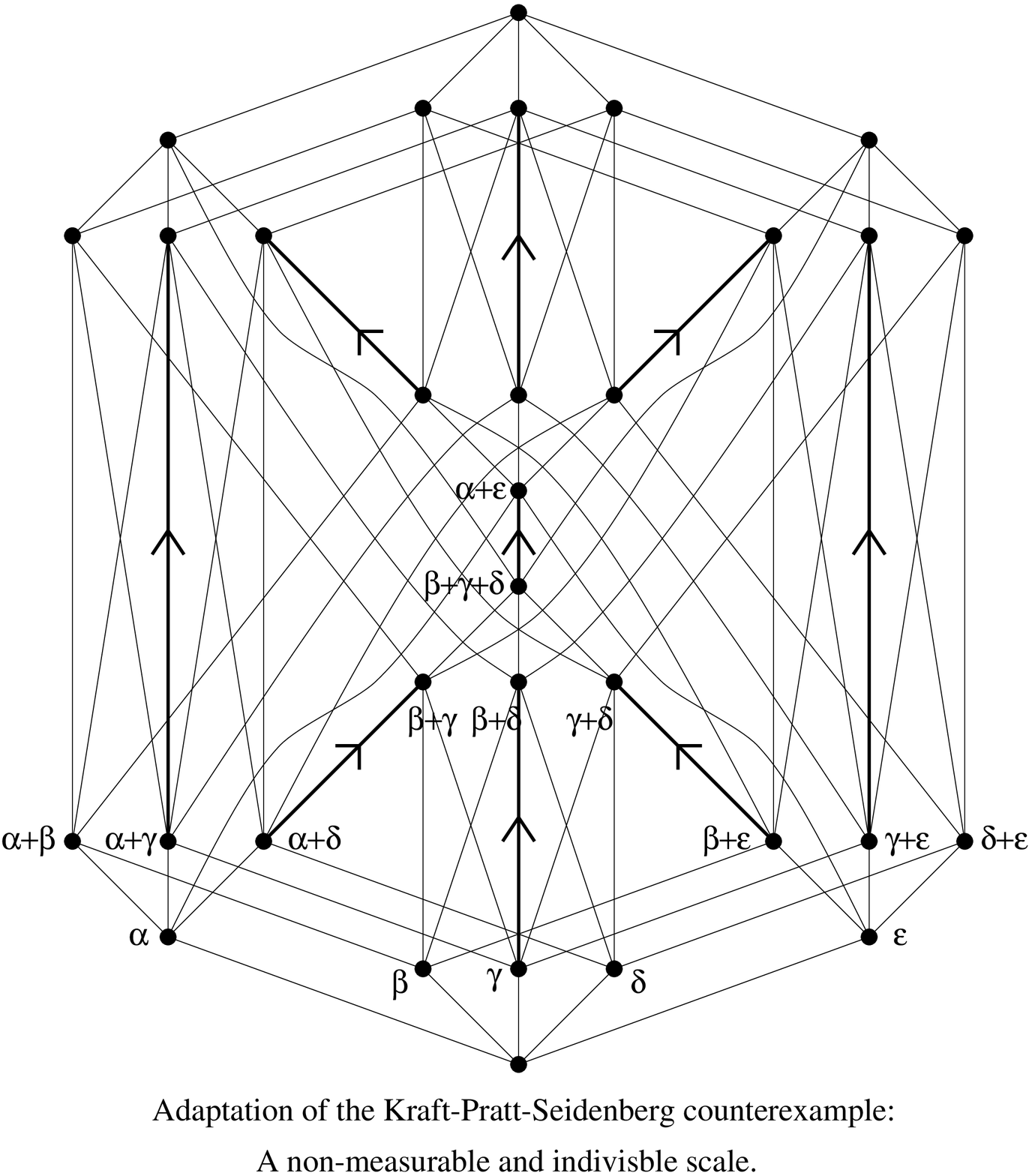}}
     \caption{}
     \label{kraftposet}
   \end{center}
 \end{figure}
my adaptation of an object constructed by Kraft, Pratt, and
Seidenberg in \cite{kraftetal}.  Their purpose was to exhibit
counterexample to Bruno de Finetti's conjecture in \cite{definetti}
that every {\em linear} ordering of a finite Boolean algebra
of propositions by comparative probabilities that satisfies the
assumptions~(\ref{prob-ordering})
of \S\ref{savages} (including weak additivity) has an
``agreeing measure.''  A broad generalization of that conjecture
states, in the language of the present paper, that every finite
scale is measurable.  This adaptation is a scale that is neither
measurable nor divisible (and so is the linearly ordered example
of which it is an adaptation).

This scale is the image of the
Boolean algebra of all subsets of $\set{a,b,c,d,e}$ under
a one-to-one scaling $\rho$ satisfying
$\rho(\set{a})=\alpha, \rho(\set{b})=\beta, \rho(\set{c})=\gamma,
\rho(\set{d})=\delta, \rho(\set{e})=\varepsilon$.  The thin lines
in Figure~\ref{kraftposet} correspond to subset relations,
e.g., a thin line goes from $\rho(\set{a})=\alpha$ to
$\rho(\set{a,e})=\alpha+\varepsilon$ because $\set{a}\subseteq\set{a,e}$.
The thick lines with arrows on them correspond to the following
inequalities not necessitated by subset relations.
\begin{equation}\label{KPSinequalities}
\begin{split}
\alpha+\delta & < \beta+\gamma \\
\beta+\varepsilon & < \gamma+\delta \\
\gamma & < \beta+\delta \\
\beta+\gamma+\delta & < \alpha+\varepsilon
\end{split}
\end{equation}
\end{example}
\begin{proposition}
Example~\ref{IndivisibleScale} is a non-measurable scale.
\end{proposition}
Although the Proposition is new, the proof is due to
Kraft, Pratt, and Seidenberg.
\begin{proof}
Suppose $\mu$ is a measure on this scale.  Then from
(\ref{KPSinequalities}) it follows that
\begin{equation}\begin{split}
\mu(\alpha)+\mu(\delta) & < \mu(\beta)+\mu(\gamma) \\
\mu(\beta)+\mu(\varepsilon) & < \mu(\gamma)+\mu(\delta) \\
\mu(\gamma) & < \mu(\beta)+\mu(\delta) \\
\mu(\beta)+\mu(\gamma)+\mu(\delta) & < \mu(\alpha)+\mu(\varepsilon).
\end{split}
\end{equation}
Unlike the addition in (\ref{KPSinequalities}),
this is old-fashioned everywhere-defined and impeccably-behaved
addition of real numbers.  Therefore we can deduce an inequality
between the sum of the four left sides and that of the
four right sides, an absurdity:
$$
\mu(\alpha)+2\mu(\beta)+2\mu(\gamma)+2\mu(\delta)+\mu(\varepsilon)
<\mu(\alpha)+2\mu(\beta)+2\mu(\gamma)+2\mu(\delta)+\mu(\varepsilon).
$$
\end{proof}
\begin{proposition}
Example~\ref{IndivisibleScale} is an indivisible scale.
\end{proposition}
\begin{proof}
Observe that all sums in this proof are sums of images
of pairwise disjoint members of a division $\mathbb{A}_1$
of $\mathbb{A}$, so by Proposition~\ref{associativity}
we have all the associativity we need.  In particular, we
will deal with two members of $\mathbb{A}_1$ whose images
under $\rho$ are equal to $\gamma$.  They are disjoint.

Assume it is divisible, so that Theorem~\ref{divided-defined}
is applicable.  Since $\gamma<\beta+\gamma+\delta<\alpha+\varepsilon
=\,\ssim(\beta+\gamma+\delta)$, it follows from
\mbox{Theorem~\ref{divided-defined} (1)} % (\ref{divided-defined1})
that the sum $\gamma+(\beta+\gamma+\delta)$ exists.
(If we didn't have Theorem~\ref{divided-defined},
we would use a more leisurely but essentially equivalent
approach:  Split the set $\set{a,e}$ according to the manner
of Example~\ref{SplittingAlgorithm} so that one of its subsets
after splitting has $\beta+\gamma+\delta$ as its image under $\rho$,
and that subset is disjoint for $\set{c}$, so the addition can
be done.  That is how Theorem~\ref{divided-defined}
saves us some work.)  Now Lemma~\ref{InequalityLemma}
can be applied, and we can add the first two of the
inequalities~(\ref{KPSinequalities}) to get
$$
\alpha+\beta+\delta+\varepsilon<\beta+\gamma+\gamma+\delta.
$$
Then \mbox{Theorem~\ref{divided-defined} (2)}
implies that we can subtract $\beta+\delta$ from both sides of this
inequality, getting 
$$
\alpha+\varepsilon<\gamma+\gamma.
$$
Since (\ref{KPSinequalities}) tells us that
$\beta+\gamma+\delta<\alpha+\varepsilon$, it
follows that
$$
\beta+\gamma+\delta<\gamma+\gamma.
$$
\mbox{Theorem~\ref{divided-defined} (2)} % (\ref{divided-defined2})
tells us we can subtract $\gamma$ from both sides of this, getting
$$
\beta+\delta<\gamma.
$$
This contradicts the inequalities (\ref{KPSinequalities}).
\end{proof}
\subsection{Divisibility and measurability}
In this section, as in \S\ref{examples} we denote the
cardinality of a set $T$ by $\left|T\right|$.
\begin{lemma}\label{cards}
Suppose $\rho:\mathbb{A}\rightarrow\mathcal{R}$ is a
divided basic scaling and $\mathbb{A}$ is the Boolean
algebra of all subsets of a finite set $S$.  Then
for $T_1,T_2\subseteq S$, if $\rho(T_1)=\rho(T_2)$
then $\left|T_1\right|=\left|T_2\right|$, and if
$\rho(T_1)<\rho(T_2)$
then $\left|T_1\right|<\left|T_2\right|$.
(The converse is false since $\rho(T_1)$ and
$\rho(T_2)$ can be incomparable.)
\end{lemma}
\begin{proof}
First assume $\rho(T_1)=\rho(T_2)$.
Proceed by induction on $\left|T_1\right|$.
If $\left|T_1\right|=0$ then the conclusion follows from
the strictly increasing nature of $\rho$.  Then suppose
$\left|T_1\right|=n+1$ and $\rho(T_1)=\rho(T_2)$.
For some $t\in T_1$, dividedness implies that
$\rho(T_1\setminus\set{t})$ must be the same as the
image under $\rho$ of some proper subset of $T_2$.
Then apply the induction hypothesis.

If $\rho(T_1)<\rho(T_2)$ then, by divideness, we can
find $T_2'\subsetneq T_2$ such that $\rho(T_1)=\rho(T_2')$.
Then proceed as above with $T_2'$ in place of $T_2$.
\end{proof}
\begin{theorem}\label{kleene}
Suppose $\rho:\mathbb{A}\rightarrow\mathcal{R}$ is a divided
basic scaling on the Boolean algebra $\mathbb{A}$ of all
subsets of a finite set $S$.  Then there is a partition
$S_1,\dots,S_k$ of $S$ such that for every
$T_1,T_2\subseteq S$, $\rho(T_1)\leq\rho(T_2)$ if and only
if $\left|T_1\cap S_i\right|\leq\left|T_2\cap S_i\right|$
for $i=1,\dots,k$.
\end{theorem}
In effect this says members of $\mathcal{R}$ can be represented
as $k$-tuples $(t_1,\dots,t_k)$, and the $i$th component $t_i$
counts the number of members of a set that are in the $i$th
equivalence class, and moreover $(t_1,\dots,t_k)\leq(u_1,\dots,u_k)$
precisely if $t_i\leq u_i$ for every $i$.  The largest member
of $\mathcal{R}$ would correspond to $(a_1,\dots,a_k)$ where each
$a_i$ is the whole number of members of the $i$th equivalence class.
\begin{proof}
For two members $s,t\in S$, Lemma~\ref{cards} implies
that $\rho(\set{s})$ and $\rho(\set{t})$ are either equal or
incomparable.  Call $s$ and $t$ equivalent iff
$\rho(\set{s})=\rho(\set{t})$, and call the equivalence
classes $S_1,\dots,S_k$.  If $T_1,T_2$ are both subsets of
the same equivalence class then by additivity
we have $\rho\left(T_1\right)=\rho\left(T_2\right)$
or $\rho\left(T_1\right)<\rho\left(T_2\right)$
according as $\left|T_1\right|=\left|T_2\right|$
or $\left|T_1\right|<\left|T_2\right|$.
More generally, additivity implies that if
$\left|T_1\cap S_i\right|\leq\left|T_2\cap S_i\right|$
for $i=1,\dots,k$ then $\rho(T_1)\leq\rho(T_2)$,
with equality between the two values of $\rho$
if and only if equality holds between the two
cardinalities for every $i\in\set{1,\dots,k}$.

Next we need to show that if for some $i,j\in\set{1,\dots,k}$
we have $\left|T_1\cap S_i\right|<\left|T_2\cap S_i\right|$
and $\left|T_2\cap S_j\right|>\left|T_1\cap S_j\right|$
then $\rho(T_1)$ and $\rho(T_2)$ are incomparable.
To see this, first create $U_1,U_2\subseteq S$ as follows.
For each $i\in\set{1,\dots,k}$ for which
$\left|T_1\cap S_i\right|<\left|T_2\cap S_i\right|$,
delete $\left|T_1\cap S_i\right|$ members from $T_2\cap S_i$,
including, but not limited to, all members of $T_1\cap T_2\cap S_i$,
to get $U_2\cap S_i$, so that $U_2$ is the union of all $k$ of these
intersections.  Similarly, for each $j\in\set{1,\dots,k}$ for which
$\left|T_1\cap S_j\right|>\left|T_2\cap S_j\right|$, delete
$\left|T_2\cap S_j\right|$ members from $T_1\cap S_j$, including,
but not limited to, all members of $T_1\cap T_2\cap S_i$,
to get $U_1\cap S_j$, so that $U_1$ is the union of all $k$ of
these intersections.  Then for each $i\in\set{1,\dots,k}$ for
which $\left|T_1\cap S_i\right|<\left|T_2\cap S_i\right|$,
we have in effect deleted {\em all} of the members of
$T_1\cap S_i$ from $T_1$, getting $U_1\cap S_1=\varnothing$,
and we have deleted the same number of members of
$T_2\cap S_i$ from $T_2$ to get $U_2\cap S_i=\varnothing$.
Since all members of $S_i$ have the same image under $\rho$,
and since in divided scales we can subtract, we have subtracted
the same thing from both sides of either the equality
$\rho(T_1)=\rho(T_2)$ or the inequality
$\rho(T_1)<\rho(T_2)$.  Therefore we must have
$\rho(U_1)=\rho(U_2)$ or $\rho(U_1)<\rho(U_2)$, according
as the equality or the inequality holds between
$\rho(T_1)$ and $\rho(T_2)$.  And there is
no $i\in\set{1,\dots,k}$ for which $S_i$ intersects
both $U_1$ and $U_2$.  If $u_1\in U_1$ then, by divisibility,
there exists $U_2'\subseteq U_2$ such that
$\rho(\set{u_1})=\rho(U_2')$.
By Lemma~\ref{cards} this implies $U_2'$ has only one
member --- call it $u_2$.  But then
$\rho(\set{u_1})=\rho(\set{u_2})$ even though
$u_1$ and $u_2$ are in different equivalence classes
-- a contradiction following from the assumption of
comparability of $\rho(T_1)$ and $\rho(T_2)$.
\end{proof}
In effect we have proved that, under the assumptions of
the theorem, $\mathcal{R}$ must be a finite ``Kleene algebra.''
This concept generalizes the concept of Boolean algebra.
A Kleene algebra is a bounded distributive lattice with a
certain sort of complementation, which is {\em not}
a ``complementation'' as the term is understood in lattice
theory.  The precise definition is:
A Kleene algebra is a partially ordered set with largest
and smallest members 1 and 0 (this is ``boundedness'')
in which any set $\set{x,y}$ of two members has an infimum
$x\vee y$ and a supremum $x\wedge y$ (i.e., it is a lattice),
and these two operations distribute over each other,
and there is an complementation $x\mapsto\ssim x$ satisfying:
\begin{eqnarray*}
\ssim 0 & = & 1, \\
\sim\ssim x & = & x, \\
\ssim(x\wedge y) & = & (\ssim x)\vee(\ssim y), \\
x\wedge\ssim x & \leq & y\vee\ssim y.
\end{eqnarray*}
In a Boolean algebra we would have $x\wedge\ssim x=0$
and $y\vee\ssim y=1$ (i.e., this would be a lattice-theoretic
complementation) instead of this weaker last condition.
A Boolean algebra can be defined as a ``complemented distributive
lattice.'' Up to isomorphism, a finite Kleene algebra is the
same thing as a family of sub-multisets of a finite multiset,
that is closed under the three operations.

{\sc Notational Convention.}  In order to use it as a tool
in the statement and proof of the next theorem, we further
develop the notation introduced in the paragraph after
Theorem~\ref{kleene}.  For any divided basic scaling
$\rho:\mathbb{A}\rightarrow\mathcal{R}$ on a finite Boolean
algebra $\mathbb{A}$, we represent members of $\mathcal{R}$
as tuples $(t_1,\dots,t_k)$ of non-negative integers.  For
any two such $k$-tuples $(t_1,\dots,t_k)$ and $(u_1,\dots,u_k)$,
we have $(t_1,\dots,t_k)\leq(u_1,\dots,u_k)$ iff
$t_i\leq u_i$ for $i=1,\dots,k$.  Addition and subtraction of
members of $\mathcal{R}$ then become term-by-term addition and
subtraction of components.  If $\alpha=(t_1,\dots,t_k)$,
$\beta=(u_1,\dots,u_k)$, $\gamma=(v_1,\dots,v_k)$, and
$\alpha\leq\beta\leq\gamma$, then the additive relative complement
$\ssim\beta_{[\alpha,\gamma]}$ is $(v_1-u_1+t_1,\dots,v_k-u_k+t_k)$.
The range $\mathcal{R}$ also has a lattice structure.  (Recall from
Example~\ref{nonarch} that a divided scale need not be a lattice if
it is not finite, and from Example~\ref{balanced} that a finite scale
need not be a lattice if it is not divided.)
The lattice structure of a finite divided scale is given
by the compontentwise definition of the meet and join operations:
\begin{eqnarray*}
(t_1,\dots,t_k)\wedge(u_1,\dots,u_k)
& = & (t_1\wedge u_1,\dots,t_k\wedge u_k) \\
(t_1,\dots,t_k)\vee(u_1,\dots,u_k)
& = & (t_1\vee u_1,\dots,t_k\vee u_k).
\end{eqnarray*}
\begin{theorem}\label{DivMeas}
Suppose $\rho:\mathbb{A}\rightarrow\mathcal{R}$ is a basic scaling
and $\mathbb{A}$ is finite.  Then $\mathcal{R}$ is divisible
if and only if it is measurable.
\end{theorem}
\begin{proof}
If $\rho_1:\mathbb{A}_1\rightarrow\mathcal{R}_1$ is a division
of $\rho:\mathbb{A}\rightarrow\mathcal{R}$, and
$\mu_1:\mathcal{R}_1\rightarrow\mathbb{R}$ is a measure,
then the restriction of $\mu_1$ to $\mathcal{R}$ is also
a measure.  Therefore no generality is lost by assuming
the scale is not just divisible, but divided, and so we
do.  Since $\mathbb{A}$ is finite, we lose no generality
by assuming $\mathbb{A}$ is the algebra of all subsets of
some finite set $S$.

Following the notation introduced in the paragraph after
Theorem~\ref{kleene}, write $\rho(T)=(t_1,\dots,t_k)$
for $T\subseteq S$.  For any $m_1,\dots,m_k>0$, the function
$\mu(T)=\sum_{i=1}^k m_i t_i$ is a measure of the sort
required, and the set of all such measures is the requisite
convex set of measures.

Conversely, assume $\mathcal{R}$ is measurable.
For every measure $\mu:\mathcal{R}\rightarrow\mathbb{R}$,
the mapping $\mu\circ\rho:\mathbb{A}\rightarrow\mathbb{R}$
is a measure on the underlying Boolean algebra.
The set of all such measures $\mu\circ\rho$ satisfying
$\mu(\rho(1))=1$ is convex and bounded.
Since it is finite, we may take $\mathbb{A}$ to be the
algebra of all subsets of a finite set $S$.
For any $\zeta,\eta\in\mathcal{R}$ satisfying $\zeta<\eta$,
and any $T_1,T_2\subseteq S$ for which $\rho(T_1)=\zeta$ and
$\rho(T_2)=\eta$ we have an inequality
\begin{equation}\label{linear1}
\sum_{t\in T_1}\mu(\rho(\set{t}))
<\sum_{t\in T_2}\mu(\rho(\set{t})).
\end{equation}
We get finitely many such inequalities, plus one equality that says
\begin{equation}\label{linear2}
\sum_{t\in S}\mu(\rho(\set{t}))=1.
\end{equation}
The solution set of the system consisting of the inequalities
(\ref{linear1}) and the equation (\ref{linear2}) in the finitely
many variables $\mu(\rho(\set{t}))$, $t\in S$, is a bounded
convex set that is the convex hull of finitely many ``corners,''
and each corner is a rational point in $\mathbb{R}^n$, where
$n=\left|S\right|$.  Let $m$ be the number of corners, and let
$M$ be the $m\times n$ matrix whose rows are the corners.
For each corner $c$ call the corresponding row of $M$ the
$c^{\text{th}}$ row, and let $d_c$ be the common denominator
of the rational numbers that are the entries in the
$c^{\text{th}}$ row.  Note that each column of $M$ corresponds
to one of the variables $\mu(\rho(\set{t}))$, and so each column
of $M$ corresponds to one of the members $t\in S$.  Call that
column the $t^{\text{th}}$ column of $M$.  Let $D$ be the diagonal
matrix whose entries are the $d_c$.  Then $M^T D$ is an integer matrix.
For each $t\in S$, let
$\sigma(\set{t})=\text{the }t^{\text{th}}\text{ row of }M^T D$.
Then, following the {\sc Notational Convention} that precedes the
statement of the Theorem, $\sigma$ is the desired divided scaling.
\end{proof}

\begin{example}
Figure~\ref{convex7} depicts the convex set of all
measures on the scale in Figure~\ref{poset6}.
   \input epsf
     \begin{figure}[htb]
     \begin{center}
     \leavevmode
     \hbox{%
     \epsfxsize=2.0in
     \epsfysize=3.3in
     \epsffile{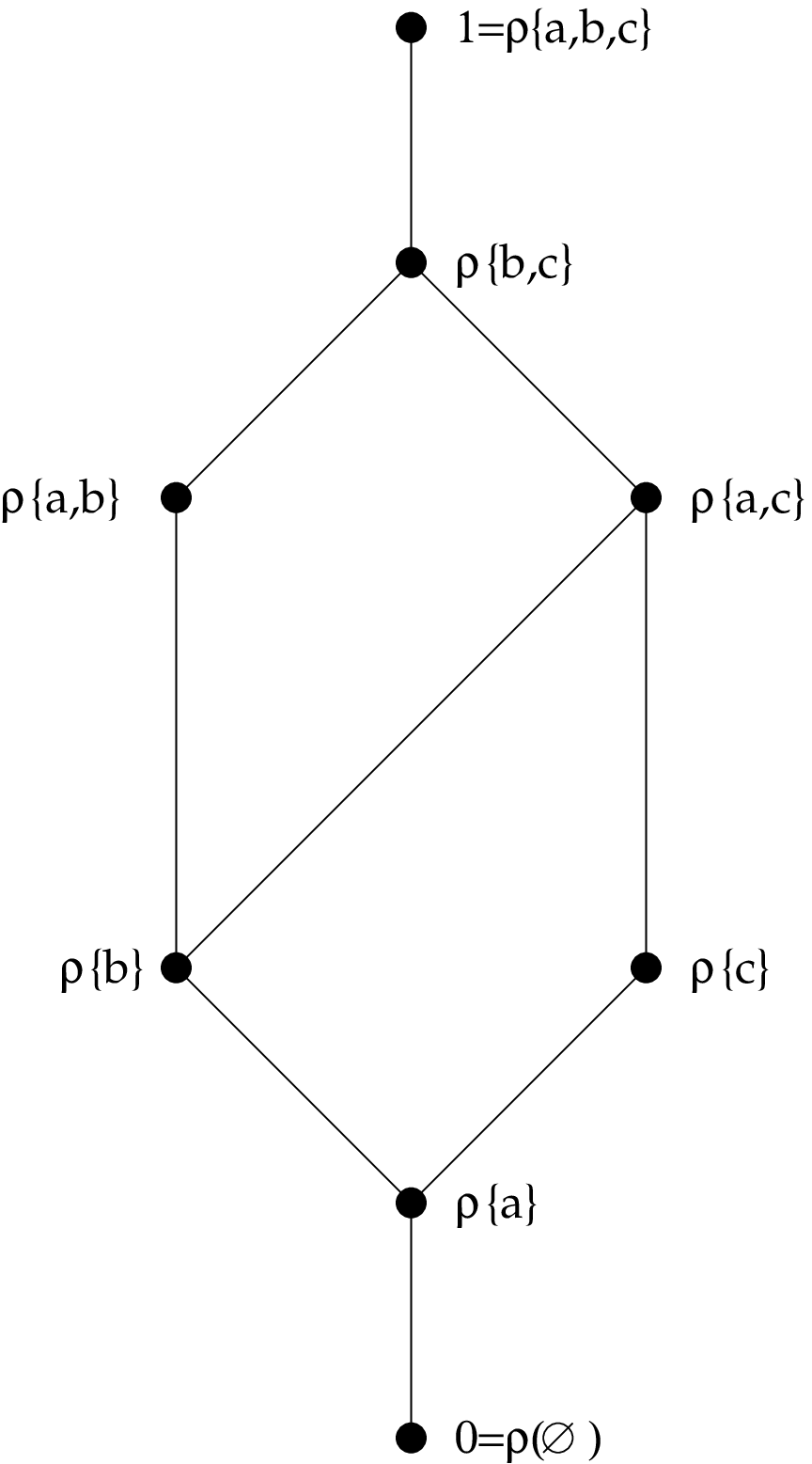}}
     \caption{}
     \label{poset6}
   \end{center}
 \end{figure}
   \input epsf
     \begin{figure}[htb]
     \begin{center}
     \leavevmode
     \hbox{%
     \epsfxsize=3.1in
     \epsfysize=2.4in
     \epsffile{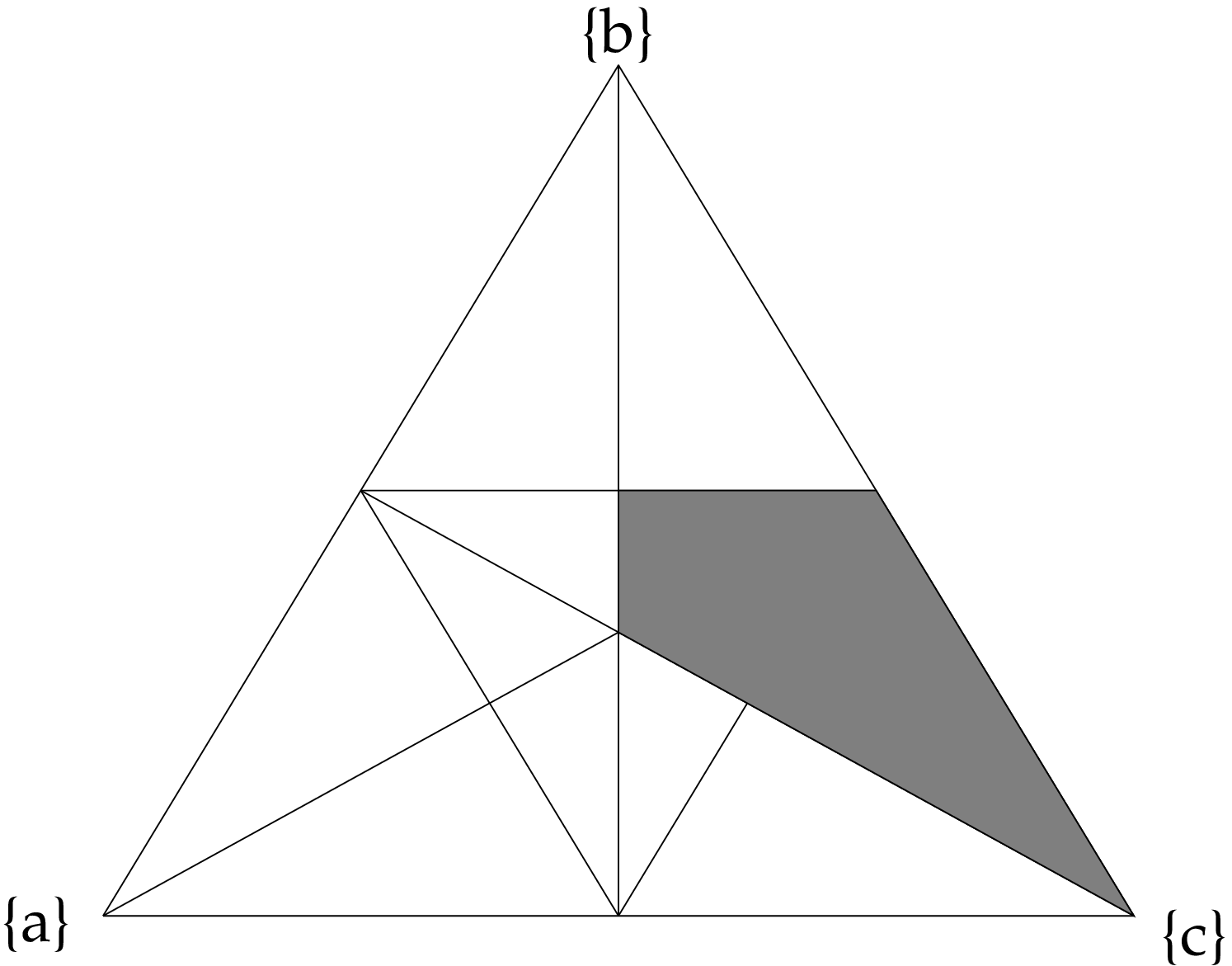}}
     \caption{}
     \label{convex7}
   \end{center}
 \end{figure}
The corners are the rows of the matrix
$$
M=\left[
\begin{array}{ccccc}
1/3 & & 1/3 & & 1/3 \\
1/4 & & 1/2 & & 1/4 \\
0   & & 1/2 & & 1/2 \\
0   & &  0  & & 1
\end{array}
\right].
$$
The respective common denominators are 3, 4, 2, and 1, and we get
\begin{eqnarray*}
\sigma(\set{a}) & = & (1,1,0,0), \\
\sigma(\set{b}) & = & (1,2,1,0), \\
\sigma(\set{c}) & = & (1,1,1,1).
\end{eqnarray*}
This means: We split the atom $a$ into two parts, and put
one in the first equivalence class and one in the second;
We split $b$ into four parts, and put one in the first equivalence
class, two in the second, and one in the third; We split $c$
into four parts, and put one in each of the four equivalence classes.
The first equivalence class has three members; the second has four;
the third has two; the fourth has one.
\end{example}

\subsection{Multiplication}\label{secmultiplication}

Theorem~\ref{divided-defined} told us that if a scale is divided,
then subtraction is generally defined, i.e., whenever
$\zeta\leq\eta$ then $\eta-\zeta$ exists.  Assume dividedness
and Archimedeanism, but replace the assumption that
$\zeta\leq\eta$, with the assumption that the scale is linearly
ordered and $\zeta\neq 0$ (so, by Archimedeanism, $\zeta$
is not infinitesimal).  Consider the maximum value of
$n\in\set{0,1,2,3,\dots}$ such that the following
difference exists
$$
\eta\underbrace{-\zeta-\zeta-\cdots-\zeta}_{n}.
$$
That maximum value may be 0, and must be finite.
We must have
$$
\alpha=\eta\underbrace{-\zeta-\zeta-\cdots-\zeta}_{n}<\zeta
$$
since otherwise we could subtract another $\zeta$, contradicting
the maximality of $n$ (here we have used the assumption of linear
ordering).  Since $\alpha<\zeta$, we can subtract $\alpha$ from
$\zeta$.  Consequently we have an algorithm:
$$
\begin{array}{lrcl}
\textsf{(1) Let} & i & = & 1\text{ (the positive integer 1,
not the maximum element of a scale)}. \\
\textsf{(2) Let} & \alpha & = & \zeta. \\
\textsf{(3) Let} & \beta & = & \eta. \\
\textsf{(4) Let} & n_i & = &
\max\left\{n\right.
:\beta\underbrace{-\alpha-\cdots-\alpha}_{
n}\left.\text{ exists}\right\} \\
\textsf{(5) Let} & \beta & = &
\beta\underbrace{-\alpha-\cdots-\alpha}_{n_i} \\
\textsf{(6) If} & \beta & = & 0\textsf{ then stop, else }
\left\{\begin{array}{l}
\textsf{Increment }i\text{ to }i+1; \\
\textsf{Interchange the values of }
\alpha\textsf{ and }\beta; \\
\textsf{Go to (4).}
\end{array}\right\}. \\
\end{array}
$$
If we had defined any reasonable notions of multiplication
and division of members of a scale, then this algorithm
would find the continued fraction expansion:
$$
\frac{\eta}{\zeta}=n_1+\frac{1}{n_2+\dfrac{1}{n_3+\dfrac{1}{n_4+
\cdots\cdots}}}
$$
Call this the formal continued fraction expansion of the
formal quotient $\eta/\zeta$.  Observe that the formal
continued fraction expansions of the formal quotients
$\eta/\zeta$ and $\theta/\zeta$ are the same only if
the difference between $\eta$ and $\theta$ is infinitesimal,
and therefore, by Archimedeanism, is 0.  All this is summarized
by a lemma:
\begin{lemma}\label{multiplicationlemma}
In a linearly ordered divided Archimedean scale, the continued
fraction and the (non-zero) value of the denominator of a formal
quotient determine the value of the numerator.
\end{lemma}
\begin{theorem}
On any linearly ordered divided Archimedean scale $\mathcal{R}$
there is exactly one measure
$\mu:\mathcal{R}\rightarrow[0,1]\subseteq\mathbb{R}$
such that $\mu(1)=1$.
\end{theorem}
\begin{proof}
Apply Lemma~\ref{multiplicationlemma} in the case
$\zeta=1\in\mathcal{R}$.  Any measure
$\mu:\mathcal{R}\rightarrow[0,1]\subseteq\mathbb{R}$
for which $\mu(1)=1$ takes addition and subtraction in
$\mathcal{R}$ to the usual addition and subtraction in
$\mathbb{R}$.  Consequently the formal continued fraction
expansion of the formal quotient $\eta/1$ must be the same
as the ordinary continued fraction expansion of $\mu(\eta)$.
The measure $\mu$ is therefore completely determined by
the structure of $\mathcal{R}$.
\end{proof}
The theorem says we can identify any divided Archimedean scale
with some subset of $[0,1]\subseteq\mathbb{R}$.
\begin{definition}\label{dfnmultiplication}
For $\zeta,\theta$ in a linearly ordered divided Archimedean
scale, the relation
$$
\zeta\theta=\eta
$$
means
$$
\mu(\zeta)\mu(\eta)=\mu(\theta),
$$
or equivalently the not-everywhere-defined multiplication
is given by
$$
\zeta\eta=\mu^{-1}(\mu(\zeta)\mu(\eta)).
$$
\end{definition}
\begin{example}
Let $\rho$ be the probability measure on the set of
all subsets of $\set{a,b,c}$ that assigns $1/3$ to
each of $\set{a},\set{b},\set{c}$.  Then
Definition~\ref{dfnmultiplication} fails to define
$\rho\set{a}\rho\set{b}$.
\end{example}

In \S\ref{productrule} we will apply
Definition~\ref{dfnmultiplication} to probability.

\section{Homomorphisms and Stone spaces}
The material in this section is not new.
All or nearly all of it can be found in \cite{halmos}.
\subsection{Homomorphisms}
\begin{definition}\label{DefnHomomorphism}
\begin{enumerate}
\item
Let $\mathbb{A},\mathbb{B}$ be Boolean algebras.
A {\bf homomorphism} $\varphi:\mathbb{A}\rightarrow\mathbb{B}$
is a function for which, for all $x,y\in\mathbb{A}$ we have:
\begin{eqnarray*}
\varphi(x\wedge y) & = & \varphi(x)\wedge\varphi(y) \\
\varphi(x\vee y) & = & \varphi(x)\vee\varphi(y) \\
\varphi(\ssim x) & = & \ssim\varphi(x).
\end{eqnarray*}
\item The {\bf kernel} of a homomorphism
$\varphi:\mathbb{A}\rightarrow\mathbb{B}$ is
$\varphi^{-1}(0)=\set{x\in\mathbb{A}:\varphi(x)=0}$.
\item A {\bf principal} homomorphism is one whose
kernel is of the form $\set{y\in\mathbb{A}:y\leq x}$
for some $x\in\mathbb{A}$.  We say that the kernel
is generated by $x$.  Other homomorphisms are
{\bf nonprincipal} homomorphisms.
\item
A homomorphism $\varphi:\mathbb{A}\rightarrow\mathbb{B}$
is {\bf 2-valued} if $\mathbb{B}$ is the two-element
Boolean algebra $\set{0,1}$.
\end{enumerate}
\end{definition}
The next proposition is an immediate corollary of
Definition~\ref{DefnHomomorphism}.
\begin{proposition}
A 2-valued homomorphism $\varphi$ on $\mathbb{A}$ is
principal if and only if for some atom $x\in\mathbb{A}$,
$\varphi(y)=$ 1 or 0 according as $x\leq y$ or $x\wedge y=0$.
\end{proposition}
\begin{example}\label{simpleexample}
Every finite Boolean algebra $\mathbb{A}$
is isomorphic to the Boolean algebra of all subsets of
some finite set $\Phi$.  Let $\mathbb{B}$ be the Boolean
algebra of all subsets of some non-empty set $\Psi\subseteq\Phi$.
For $x\in\mathbb{A}$, let $\varphi(x)=x\cap\Psi$.  Then
$\varphi$ is a principal homomorphism whose kernel is
generated by $\Phi\diagdown\Psi$.  If $\Psi$ is a single-element
set, then $\varphi$ is a 2-valued homomorphism.
\end{example}
\begin{example}  Let $\mathbb{A}$ be the Boolean algebra
of all subsets of $\mathbb{N}=\set{0,1,2,3,\dots}$.
Let $\mathbb{B}$ be the Boolean algebra of all equivalence
classes of such sets, two sets $A,B\subseteq\mathbb{N}$
being considered equivalent if
$\left|(A\diagdown B)\cup(B\diagdown A)\right|<\aleph_0$
(this is a much coarser equivalence relation than the one
considered in Example~\ref{nonarch}!).  Meet, join, and complement
on $\mathbb{B}$ are defined by choosing members of equivalence
classes, then evaluating the meet, join, or complement of
those, then taking the equivalence class to which the result
belongs.  It is easy to check that these operations are
well-defined.  For $x\in\mathbb{A}$, let $\varphi(x)$ be
the equivalence class to which $x$ belongs.  This is a
nonprincipal homomorphism whose kernel is the set of all
finite subsets of $\mathbb{N}$.
\end{example}
\begin{example} Let $\mathbb{A}$ be the Boolean algebra of
all subsets of $\mathbb{N}$ that are either finite or cofinite
(cofinite means having a finite complement).  For any
$n\in\mathbb{N}$ the mapping
$$
x\mapsto\left\{\begin{array}{cl}
1 & \text{if }n\in x \\
0 & \text{if }n\not\in x
\end{array}\right\}
$$
is a principal 2-valued homomorphism whose kernel
is the set of all subsets of $\mathbb{N}\diagdown\set{n}$.
The mapping
$$
x\mapsto\left\{\begin{array}{cl}
1 & \text{if }x\text{ is cofinite} \\
0 & \text{if }x\text{ is finite}
\end{array}\right\}
$$
is a nonprincipal 2-valued homomorphism whose kernel
is the set of all finite subsets of $\mathbb{N}$.
\end{example}
\begin{example}\label{topologicalexample}
Let $\mathbb{A}$ be the Boolean algebra
of all clopen (i.e., simultaneously closed and open) subsets
of the Cantor set $C$.  Let $\mathbb{B}$ be the Boolean algebra %(
of all clopen subsets of $C\diagdown[0,1/3)$. %](
The mapping $x\mapsto x\diagdown[0,1/3)$ %]
is a nonprincipal homomorphism whose kernel is %(
$\set{x\in\mathbb{A}:x\subseteq[0,1/3)}$. %]
Now fix one point $r\in C$.
The mapping
\begin{equation}\label{typicalstonepoint}
x\mapsto\left\{\begin{array}{cl}
1 & \text{if }r\in x \\
0 & \text{if }r\not\in x
\end{array}\right\}
\end{equation}
is a nonprincipal 2-valued homomorphism whose
kernel is $\set{x\in\mathbb{A}:r\not\in x}$.
\end{example}

\subsection{Stone's representation of Boolean algebras}

In both Example~\ref{simpleexample}
and Example~\ref{topologicalexample}, we saw a 2-valued
homomorphism on a Boolean algebra $\mathbb{A}$ of sets defined
as in (\ref{typicalstonepoint}) -- its value is 1 or 0 according
as the set does or does not contain a certain point.
In Example~\ref{simpleexample}, the homomorphism was principal
because the set containing only that one point was a member of
$\mathbb{A}$; in Example~\ref{topologicalexample}, it was
nonprincipal because the set containing only that one point
was not a member of $\mathbb{A}$.  We shall see that in a sense,
these examples are typical of 2-valued homomorphisms: We can
represent an arbitrary Boolean algebra $\mathbb{A}$ as the Boolean
algebra of certain subsets of a certain set $\Phi(\mathbb{A})$,
and then find that every 2-valued homomorphism is of the form
(\ref{typicalstonepoint}).
The homomorphism will be principal or nonprincipal according
as the set containing only the point that so represents it is
or is not one of the subsets of $\Phi(\mathbb{A})$ that are
identified with members of the Boolean algebra $\mathbb{A}$.

So we let
$$
\Phi(\mathbb{A})
=\text{the set of all 2-valued homomorphisms on }\mathbb{A},
$$
and we identify each $x\in\mathbb{A}$ with
\begin{equation}\label{typicalclopenset}
\set{\varphi\in\Phi(\mathbb{A}):\varphi(x)=1}
\end{equation}
$=$ the set of all 2-valued homomorphisms on $\mathbb{A}$
that map $x$ to 1.  The operations of meet, join, and complement
in $\mathbb{A}$ correspond to the operations of (finite)
intersection, (finite) union, and set-theoretic
complementation on subsets of $\Phi(\mathbb{A})$.
(Infinitary operations are more problematic.  The infinitary
join of $\mathbb{X}\subseteq\mathbb{A}$ is
the smallest upper bound $\bigvee\mathbb{X}$ of
$\mathbb{X}$ in $\mathbb{A}$.  This join does not always
exist -- counterexamples can be found within
Example~\ref{topologicalexample}.  When the join does exist,
it does not generally correspond to the union of 
$\bigcup_{x\in\mathbb{X}}\set{\varphi\in\Phi(\mathbb{A}):\varphi(x)=1}$,
since the union of sets of the form (\ref{typicalclopenset})
is not generally of the form (\ref{typicalclopenset}).
Rather, the join corresponds to the smallest set of the
form (\ref{typicalclopenset}) that includes the union.)

The next result is {\em Stone's representation theorem}.
\begin{theorem}
The mapping
$$
x\mapsto\Phi(x)=\set{\varphi\in\Phi(\mathbb{A}):\varphi(x)=1}
$$
is an isomorphism from the Boolean algebra $\mathbb{A}$,
to the Boolean algebra of sets of the form (\ref{typicalclopenset})
with the operations of intersection, union, and set-theoretic
complementation in the roles of meet, join, and complement.
\end{theorem}
\begin{proof}
First we show that $\Phi$ is a homomorphism.
\begin{eqnarray*}
\Phi(\ssim x) & = &
\set{\varphi\in\Phi(\mathbb{A}):\varphi(\ssim x)=1} \\
& = & \set{\varphi\in\Phi(\mathbb{A}):\varphi(x)\neq 1}
\text{since }\varphi\text{ is a 2-valued homomorphism,} \\
& = & \Phi(\mathbb{A})\diagdown\set{\varphi\in\Phi:\varphi(x)=1} \\
& = & \Phi(\mathbb{A})\diagdown\Phi(x).  \\
& & \text{So complements in }\mathbb{A}
\text{ go to set-theoretic complements.}  \\
\Phi(x\wedge y) & = &
\set{\varphi\in\Phi(\mathbb{A}):\varphi(x\wedge y)=1} \\
& = & \set{\varphi\in\Phi(\mathbb{A}):\varphi(x)\wedge\varphi(y)=1}
\text{since }\varphi\text{ is a homomorphism,} \\
& = & \set{\varphi\in\Phi(\mathbb{A}):\varphi(x)=1\text{ and }
\varphi(y)=1} \\
& = & \set{\varphi\in\Phi(\mathbb{A}):\varphi(x)=1}
\cap\set{\varphi\in\Phi(\mathbb{A}):\varphi(y)=1} \\
& = & \Phi(x)\cap\Phi(y). \\
& & \text{So meets in }\mathbb{A}\text{ go to intersections.}
\end{eqnarray*}
Let ``$\vee$'', ``or'', and ``$\cup$'' replace
``$\wedge$'', ``and'', and ``$\cap$'' respectively, to show that
joins in $\mathbb{A}$ go to unions.
So $\Phi$ is indeed a homomorphism.

To show that it is an isomorphism, we need to show
that it is one-to-one.  For $x,y\in\mathbb{A}$ let
$$
x+y=(x\wedge\ssim y)\vee(y\wedge\ssim x)
=(x\vee y)\wedge\ssim(x\wedge y),
$$
and let
$$
xy=x\wedge y.
$$
Then it can be checked that $\mathbb{A}$ becomes a
commutative ring with zero element 0 and unit element 1,
in which every element is idempotent and every element
is its own additive inverse.  The Boolean operations of
meet, join, and complement can be recovered from the
ring operations:
\begin{eqnarray*}
x\wedge y & = & xy, \\
x\vee y & = & x+y+xy, \\
\ssim x & = & 1+x.
\end{eqnarray*}
And Boolean homomorphisms coincide exactly with ring homomorphisms.
The kernel $\set{x\in\mathbb{A}:\Phi(x)=\varnothing}$
of the Boolean homomorphism is the same thing as the kernel of
the ring homomorphism.  Therefore, to show that $\Phi$ is one-to-one,
it is enough to show that the kernel contains only $0\in\mathbb{A}$.
That is the same as showing that if $x>0$ then
$\Phi(x)\neq\varnothing$.  In other words, if $x\neq 0$ then
for some 2-valued homomorphism $\varphi$ on $\mathbb{A}$
we have $\varphi(x)=1$.  Equivalently, if $x\neq 1$ then
for some 2-valued homomorphism $\varphi$ on $\mathbb{A}$
we have $\varphi(x)=0$.  The kernel of such a homomorphism
is a proper ideal.  The fact that it has only two cosets
implies that it is a {\em maximal} proper ideal.  So we
need only show that any $x\neq 1$ is a member of some
maximal proper ideal.  That is well-known to follow from
a standard application of Zorn's lemma.
\end{proof}

\subsection{Topology}

As in Example~\ref{topologicalexample} these sets
$\set{\varphi\in\Phi:\varphi(x)=1}$
will be the clopen subsets of $\Phi$ -- but to say that,
we need a topology on $\Phi$.  Here it is.
\begin{definition}
The {\bf Stone space} $\Phi(\mathbb{A})$ of a Boolean
algebra $\mathbb{A}$ is the set $\Phi$ of all 2-valued
homomorphisms on $\mathbb{A}$ endowed with the topology
whose basic open sets are sets of the form (\ref{typicalclopenset}).
That means the open sets are just those that are unions
of arbitrary collections of sets of the form (\ref{typicalclopenset}).
\end{definition}
This is the same as the topology of pointwise convergence of
nets of homomorphisms.  That the basic open sets are closed
follows immediately from the fact that the basic open set
$\Phi(\ssim x)$ is complementary to the basic open set $\Phi(x)$.
That the basic open sets are the {\em only} clopen sets is
proved in \cite{halmos} by using the following theorem.
But our real motive for including this theorem is its
use in \S\S\ref{continuity} and \ref{bayes}.
\begin{theorem}\label{compactness}
(Compactness) 
Let $\mathbb{X}\subseteq\mathbb{A}$. Suppose for every finite
subset $\mathbb{X}_0\subseteq\mathbb{X}$ there is a 2-valued
homomorphism $\varphi$ on $\mathbb{A}$ such
that for every $x\in\mathbb{X}_0$, $\varphi(x)=1$.  Then there
is a 2-valued homomorphism $\varphi$ on $\mathbb{A}$ such that
for every $x\in\mathbb{X}$, $\varphi(x)=1$.

In other words, the Stone space is compact.
\end{theorem}
\begin{proof}  We follow closely the argument in
\cite{halmos} pp.~77-78.  It suffices to prove
$\Phi(\mathbb{A})$ is a closed subset of the space $\Omega$
of {\em all} functions (not just homomorphisms) from
$\mathbb{A}$ into $\set{0,1}$, with the product topology,
since that is a compact Hausdorff space.  We have
\begin{eqnarray*}
\Phi(\mathbb{A}) & = & \left(\bigcap_{x\in\mathbb{A}}
\set{\varphi\in\Omega: \varphi(\ssim x)=\ssim\varphi(x)}\right) \\
     & \cap  & \left(\bigcap_{x,y\in\mathbb{A}}
\set{\varphi\in\Omega: \varphi(x \vee y)
=\varphi(x)\vee\varphi(y)}\right) \\
     & \cap & \left(\bigcap_{x,y\in\mathbb{A}}
\set{\varphi\in\Omega: \varphi(x \wedge y)
=\varphi(x)\wedge\varphi(y)}\right).
\end{eqnarray*}
This is closed if the sets whose intersection is taken
are closed.  They are closed because $\varphi(x)$
depends continuously on $\varphi$.
\end{proof}
\begin{example}
If $\mathbb{A}$ is the Boolean algebra of all finite
or cofinite subsets of $\mathbb{N}$, then
$\Phi(\mathbb{A})$ is the one-point compactification of
the discrete space whose underlying set is $\mathbb{N}$.
The isolated points of $\Phi(\mathbb{A})$ correspond to
principal 2-valued homomorphisms.  The one limit point
corresponds to the one nonprincipal 2-valued homomorphism,
which maps cofinite sets to 1 and finite sets to 0.
\end{example}
\begin{example}
If $\mathbb{A}$ is the Boolean algebra of all subsets of
$\mathbb{N}$, then $\Phi(\mathbb{A})$ is the Stone-Cech
compactification of the discrete space whose underlying
set is $\mathbb{N}$.  Again, the isolated points correspond
to the principal 2-valued homomorphisms, and the $2^{2^{\aleph_0}}$
limit points to the nonprincipal 2-valued homomorphisms.
If $\mathbb{B}$ is the Boolean algebra of equivalence classes
of such sets, where two sets are equivalent if and only if
their symmetric difference is finite, then the Stone space
$\Phi(\mathbb{B})$ of this atomless Boolean algebra
is the set of all limit points of $\Phi(\mathbb{A})$.
\end{example}
\begin{example}
If $\mathbb{A}$ is, as in Example~\ref{topologicalexample},
the Boolean algebra of all clopen subsets of the Cantor set $C$,
then $\Phi(\mathbb{A})=C$.
\end{example}

\section{Continuity of scalings}\label{continuity}

\subsection{Definition and examples}

\begin{definition}
A scaling $\rho:\mathbb{A}\rightarrow\mathcal{R}$
is {\bf continuous} at a homomorphism
$\varphi:\mathbb{A}\rightarrow\mathbb{B}$ if
\begin{equation}
\rho\left(\bigwedge\set{x:\varphi(x)=1}\right)
=\bigwedge\set{\rho(x):\varphi(x)=1}
\end{equation}
or, equivalently
\begin{equation}\label{second}
\rho\left(\bigvee\set{x:\varphi(x)=0}\right)
=\bigvee\set{\rho(x):\varphi(x)=0}.
\end{equation}
\end{definition}
We shall see that continuity at every principal homomorphism
is like ``continuity of measure,'' and continuity at every
2-valued homomorphism at least sometimes entails ``Archimedeanism.''
\begin{example}
If $\mathbb{A}$ is finite, then every scaling on $\mathbb{A}$
is continuous at every homomorphism on $\mathbb{A}$.
\end{example}
\begin{example}
Let $\varphi$ be the canonical homomorphism from the Boolean
algebra of all subsets of $\mathbb{N}$ into quotient algebra
of that Boolean algebra by the ideal of finite subsets of
$\mathbb{N}$.  In other words, for $A\subseteq\mathbb{N}$
we have $\varphi(A)=0$ if and only if $A$ is finite, or,
equivalently, for $A, B\subseteq\mathbb{N}$ we have
$\varphi(A)=\varphi(B)$ if and only if the symmetric
difference $(A\diagdown B)\cup(B\diagdown A)$ is finite.
Let $\rho$ be the ``simple non-Archimedean scaling'' of
Example~\ref{nonarch}.  Then $\rho$ is discontinuous at $\varphi$.
To see this, observe that
\begin{eqnarray*}
\rho\left(\bigwedge
\set{\set{n, n+1, n+2,\dots}: n\in\mathbb{N}}\right)
=\rho(\varnothing) \\
<\rho(\varnothing)+5\leq \rho(\set{n, n+1, n+2,\dots})
\quad\text{ for every }n\in\mathbb{N}.
\end{eqnarray*}
\end{example}
\begin{example}\label{nearlebesgue} (A completely additive measure)
Let $\mathbb{A}$ be the quotient algebra of Lebesgue-measurable
subsets of the interval $[0,1]$ on the real line by the ideal
of sets of measure 0.  Let $\rho$ be the quotient measure of
Lebesgue measure on $\mathbb{A}$.  This scaling is continuous
at all principal homomorphisms on $\mathbb{A}$.
\end{example}
\subsection{Continuity and additivity}
\begin{definition}
A subset $\mathbb{X}$ of a Boolean algebra $\mathbb{A}$ is
{\bf pairwise disjoint} if any distinct $x,y\in\mathbb{X}$
are disjoint, i.e., for any $x,y\in\mathbb{X}$, if $x\neq y$
then $x\wedge y=0$.
\end{definition}
\begin{definition}
Suppose $\mathbb{X}\subseteq\mathbb{A}$ is pairwise
disjoint.  Then the sum on the left side of the
equality below is defined to be the join on the
right.  The sum exists whenever the join exists.
$$\sum_{x\in\mathbb{X}}\rho(x)
=\bigvee\set{\sum_{x\in\mathbb{X}_0}\rho(x)
:\mathbb{X}_0\finsub\mathbb{X}}.$$
\end{definition}
\begin{definition}
A scaling $\rho:\mathbb{A}\rightarrow\mathcal{R}$ is
{\bf completely additive} if for every pairwise disjoint
$\mathbb{X}\subseteq\mathbb{A}$ possessing a join
$\bigvee\mathbb{X}\in\mathbb{A}$, we have
$$
\sum_{x\in\mathbb{X}}\rho(x)
=\rho\left(\bigvee\mathbb{X}\right).
$$
\end{definition}
The scaling of Example~\ref{nearlebesgue} is a completely
additive measure.  Lebesgue measure itself is only countably,
and not completely, additive.  In this example, sets of
measure zero all belong to the same equivalence class, which
is the zero-element of the quotient algebra.  Consequently we
cannot have any uncountable antichain (an ``antichain'' is
pairwise disjoint collection of members of a poset) whose
join is 1.  Only such a collection could serve as the needed
counterexample to complete additivity.
\begin{example}\label{identity}
(A completely additive scaling on a Boolean
algebra that does not satisfy the countable antichain condition)
Let $\rho$ be the identity mapping on the Boolean algebra of all
subsets of the real line.  Clearly $\rho$ is completely additive.
This Boolean algebra has uncountable antichains, i.e., it
does not satisfy the ``countable antichain condition.''
\end{example}
\begin{theorem}
Suppose a scaling $\rho:\mathbb{A}\rightarrow\mathcal{R}$
is continuous at every principal homomorphism on $\mathbb{A}$.
Then $\rho$ is completely additive.
\end{theorem}
\begin{proof}
The problem is to show that if $\mathbb{X}\subseteq\mathbb{A}$
is pairwise disjoint and has a join in $\mathbb{A}$ then
$$
\rho\left(\bigvee\mathbb{X}\right)
=\sum_{x\in\mathbb{X}}\rho(x)
=\bigvee\set{\sum_{x\in\mathbb{X}_0}\rho(x)
:\mathbb{X}_0\finsub\mathbb{X}}.
$$
Observe that
$$
\bigvee\mathbb{X}
=\bigvee\set{\bigvee\mathbb{X}_0:\mathbb{X}_0\finsub\mathbb{X}}.
$$
Therefore by condition (\ref{second}), characterizing continuity,
and the assumption that $\rho$ is continuous at every principal
homomorphism, it suffices that there be a principal homomorphism
whose kernel is $[0,\bigvee\mathbb{X}]$.  That homomorphism is
$x\mapsto x\wedge\ssim\bigvee\mathbb{X}$
from $\mathbb{A}$ into the relative Boolean algebra
$\left[0, \ssim\bigvee\mathbb{X}\right]$.
\end{proof}
\subsection{Continuity and Archimedeanism}
\begin{definition}\label{DefInfinitesimal}
A member $\delta\in\mathcal{R}$ is an {\bf infinitesimal}
for a basic scaling $\rho:\mathbb{A}\rightarrow\mathcal{R}$
if for some infinite pairwise disjoint
$\mathbb{X}\subseteq\mathbb{A}$ we have
$\rho(x)\geq\delta$ for every $x\in\mathbb{X}$.
\end{definition}
\begin{example}
The zero element of any scale is an infinitesimal.
\end{example}
\begin{example}
If $\rho:\mathbb{A}\rightarrow\mathcal{R}=[0,1]\subseteq\mathbb{R}$
is a measure, then there are no nonzero infinitesimals
in $\mathcal{R}$.
\end{example}
\begin{example} In Example~\ref{nonarch},
every member of the ``initial galaxy''
$$
\set{\rho(\varnothing),\,\rho(\varnothing)+1,\,
\rho(\varnothing)+2,\,\rho(\varnothing)+3,\,\dots\dots}
$$
is an infinitesimal.
\end{example}
\begin{example}
The identity mapping from any Boolean algebra to itself
is a basic scaling; the algebra regarded as a scale has no
infinitesimals.
\end{example}
\begin{proposition}
Suppose $\rho:\mathbb{A}\rightarrow\mathcal{R}$ is a basic 
scaling and $\sigma:\mathcal{R}\rightarrow\mathcal{S}$ is
a scaling.  (Recall that according to
Definition~\ref{extensiondefinition}, $\sigma$ ``extends''
$\rho$.)  If $\delta\in\mathcal{R}$ is an infinitesimal,
then so is $\sigma(\delta)\in\mathcal{S}$.
\end{proposition}
It is easy to see that the converse is false:
\begin{example}
The Boolean algebra of all subsets of any set, viewed as a
scale, contains no infinitesimals.
\end{example}
In other words, extending a scale can create infinitesimals
but cannot destroy them.

So now we have motivated the next definition.
\begin{definition}\label{DefArchimedean}
\begin{enumerate}
\item
Let $\rho:\mathbb{A}\rightarrow\mathcal{R}$ be a divided
basic scaling.  The scale $\mathcal{R}=\set{\rho(x):x\in\mathbb{A}}$
is {\bf Archimedean} if it has no nonzero infinitesimals,
and {\bf non-Archimedean} if it contains at least one
nonzero infinitesimal.
\item An Archimedean divided scale is {\bf stably Archimedean}
if there is no scaling $\sigma:\mathcal{R}\rightarrow\mathcal{S}$
extending the basic scaling $\rho$, such that $\mathcal{S}$
contains any nonzero infinitesimal, and {\bf unstably Archimedean}
if it is Archimedean but not stably Archimedean.
\end{enumerate}
\end{definition}
%
%  Tim Chow suggested "stably Archimedean"       3/26/01
%  Dan Luecking suggested "stably Archimedean"   3/27/01
%
The term ``stably Archimedean'' was suggested by Timothy Chow
and Daniel Lueking independently of each other, in response to
a request for suggested nomenclature posted to the usenet
newsgroup \textsf{sci.math.research}.

Why does Definition~\ref{DefArchimedean} say ``divided''?
Suppose $\rho(x)=\alpha<\beta=\rho(y)$, and
there is some infinite pairwise disjoint collection
$\mathbb{U}$ of members of the domain of $\rho$
such that for any $u\in\mathbb{U}$ we have
$u\wedge x=0$ and $\rho(x)+\rho(u)\geq\rho(y)$.
Divideness implies we can subtract $\alpha$ from
$\beta$, and Definition~\ref{DefInfinitesimal}
then implies $\beta-\alpha$ is an infinitesimal.
Without dividedness I see no way to guarantee that
any nonzero lower bound of $\set{\rho(u):u\in\mathbb{U}}$
exists.  Thus, without divideness, it is conceivable
that two members of a scale could differ infinitesimally,
even though no member differs infinitesimally from 0.
I am indebted to an anonymous referee for this point.
If I knew any such example, I would consider emending
Definition~\ref{DefArchimedean}.
(The referee speculated that if $F$ is an ordered field of
which the real field $\mathbb{R}$ is a subfield, so that
$F$ contains infinitesimals, then the
set
$$
\set{\alpha\in F: 0\leq\alpha\leq 1
\text{ and neither }\alpha\text{ nor }1-\alpha
\text{ is a nonzero infinitesimal}}
$$
would be such a case.  But it must be remembered that, by
our definitions, the addition on a scale is inherited from
a scaling whose domain is some Boolean algebra.  No such
mapping was proposed.)

\begin{theorem}\label{archimedean}
Suppose a scaling $\rho:\mathbb{A}\rightarrow\mathcal{R}$
is continuous at every 2-valued homomorphism on $\mathbb{A}$,
and $\mathcal{R}$ is linearly ordered.  Then $\mathcal{R}$
has no infinitesimals.  (Consequently, if $\mathcal{R}$ is
divided, it is Archimedean.)
\end{theorem}
\begin{proof}
Suppose $\delta>0$ is an infinitesimal in $\mathcal{R}$.
We have seen % in the section on Stone spaces
that for any 2-valued homomorphism $\varphi$,
the infimum $\bigwedge\set{x:\varphi(x)=1}$ exists, and
$$
\bigwedge\set{x:\varphi(x)=1}
=\left\{
\begin{array}{cl}
\text{an atom }x_\varphi & \text{if }\varphi
\text{ is principal} \\
0 & \text{if }\varphi\text{ is nonprincipal}
\end{array}
\right\}.
$$
In the principal case, for each $x\in\mathbb{A}$
we have
$$
\varphi(x)=
\left\{
\begin{array}{cl}
1 & \text{if }x_\varphi\leq x \\
0 & \text{if }x_\varphi\not\leq x
\end{array}
\right\}.
$$
We have now defined $x_\varphi$ when $\varphi$ is
a principal homomorphism; next we shall define
$x_\varphi$ in terms of $\delta$ when $\varphi$
is a nonprincipal homomorphism.  In the latter case,
since the greatest lower bound
$\bigwedge\set{\rho(x):\varphi(x)=1}$ is $0<\delta$,
it must be that $\delta$ is not a lower bound,
and that means some $x\in\mathbb{A}$ satisfies
$\varphi(x)=1$ and $\rho(x)<\delta$.
Choose such an $x$ and call it $x_\varphi$.
Via Stone's duality we can identify $x_\varphi$ with
a clopen subset of the Stone space --- the set of all
2-valued homomorphisms that map $x_\varphi$ to 1 ---
which contains the point $\varphi$.  Now we have a
clopen cover $\set{x_\varphi:\varphi\in\Phi}$ of the
Stone space.  Since the Stone space is compact, this
has a finite subset $\set{x_{\varphi_1},\dots\dots,x_{\varphi_n}}$
that covers the whole Stone space, so that
$x_{\varphi_1}\vee\cdots\cdots\vee x_{\varphi_n}=1$.
Some terms in this join -- call them
$x_{\varphi_{m+1}},\dots\dots,x_{\varphi_n}$
-- may be atoms whose images under $\rho$ are $\geq\delta$.
The join $x_{\varphi_1}\vee\cdots\cdots\vee x_{\varphi_m}$
of the others must be $\geq$ any $x\in\mathbb{A}$ whose
image under $\rho$ is $<\delta$.  Since these $x$'s
need not be disjoint, we replace them with $y_1,\dots\dots,y_m$
such that $y_i\leq x_{\varphi_i}$ for $i=1,\dots\dots m$,
$y_i\wedge y_j=0$ for $i,j=1,\dots\dots m$, and
$y_1\vee\cdots\cdots\vee y_m
=x_{\varphi_1}\vee\cdots\cdots\vee x_{\varphi_m}$.
(This can be done by letting
$y_i=x_{\varphi_i}\wedge\ssim(\cdots\cdots\vee x_{\varphi_{i-1}})$
for each $i$.)

That $\delta$ is an infinitesimal means there is an infinite
pairwise disjoint set $\mathbb{Z}\subseteq\mathbb{A}$
such that for each $z\in\mathbb{Z}$ we have $\rho(z)\geq\delta$.
Since the complement of $y_1\vee\cdots\cdots\vee y_m$ consists
of only finitely many atoms, no generality is lost by assuming
$$
\bigvee_{z\in\mathbb{Z}}z\leq y_1\vee\cdots\cdots\vee y_m.
$$
This inequality entails
$$
\rho\left(\bigvee_{z\in\mathbb{Z}}z\right)
\leq\rho\left(y_1\vee\cdots\cdots\vee y_m\right),
$$
and that in turn entails the middle inequality below:
$$
\sum_{z\in\mathbb{Z}}\delta
\leq\sum_{z\in\mathbb{Z}}\rho(z)
\leq\sum_{i=1}^m\rho(y_i)\leq
\underbrace{\delta+\cdots\cdots+\delta}_{m\text{ terms}}.
$$
Since the first sum has infinitely many terms and $\delta>0$,
that is not consistent with Lemma~\ref{InequalityLemma}
(\ref{InequalityLemma1}).
\end{proof}
\begin{theorem}\label{measurable}
If a scale is linearly ordered, Archimedean,
and divided, then it is measurable.
\end{theorem}
\begin{proof}
For any $\alpha,\beta\in\mathcal{R}$, linear ordering implies
that either $\alpha<\ssim\beta$, $\alpha>\ssim\beta$, or
$\alpha=\ssim\beta$.  Dividedness then entails that in the
first case, $\alpha+\beta$ exists and $\alpha\oplus\beta$
does not, in the second case $\alpha\oplus\beta$ exists and
$\alpha+\beta$ does not, and in the third case they both exist,
and $\alpha+\beta=1$ and $\alpha\oplus\beta=0$.

We define an abelian group $G$ whose underlying set
is $\mathbb{Z}\times\mathcal{R}$, i.e., the set of all
ordered pairs $(n,\alpha)$ where $n$ is an integer and
$\alpha\in\mathcal{R}$, modulo the identification of
$(n+1,0)$ with $(n,1)$, for each $n\in\mathbb{Z}$.
The addition in this group is
\begin{equation}
(n,\alpha)+(m,\beta)=\left\{
\begin{array}{ll}
(n+m,\alpha+\beta) & \text{if }\alpha\leq\ssim\beta \\
(n+m+1,\alpha\oplus\beta) & \text{if }\alpha\geq\ssim\beta
\end{array}
\right.
\end{equation}
The identification of $(n+1,0)$ with $(n,1)$ keeps
the two pieces of this definition from contradicting
each other.  We linearly order this group by saying
that if $\alpha,\beta\neq 1$ then $(n,\alpha)<(m,\beta)$
if either $n<m$, or $n=m$ and $\alpha<\beta$.
This linear ordering is compatible with the group
operation, in the sense that for any $0\neq u\in G$, either
$u>0$ or $-u>0$, and for any $u,v,w\in G$, if $u<v$ then
$u+w<v+w$.  A group with such a compatible linear ordering
is a ``linearly ordered group.''

Observe that the Archimedean nature of $\mathcal{R}$ and
that of $\mathbb{Z}$ together imply that $G$ is Archimedean
in the sense that for any $u,v>0$ in $G$, there is some
positive integer $n$ such that
$$
\underbrace{u+\cdots\cdots+u}_{n\text{ terms}}>v,
$$
so that no matter how small $u$ is by comparison to $v$,
it takes only finitely many $u$'s to add up to more than $v$.

A well-known theorem of H\"older (see \cite{fuchs}, p.~45)
says that if a linearly ordered group $G$ is Archimedean,
then there is an isomorphism $f$ from $G$ into the additive
group of real numbers.  For $\alpha\in\mathcal{R}$, so that
$(0,\alpha)\in G$, let $\mu(\alpha)=f(0,\alpha)$.
Then $\mu$ is the desired measure.
\end{proof}
Note that any extension of Example \ref{IndivisibleScale}
to a linearly ordered scale is a counterexample showing that
the hypothesis of divisibility cannot be dispensed with.

I do not know how to prove the following.

{\sc Conjecture}. The hypothesis of linear ordering
in Theorem~\ref{measurable} can be dropped.

\section{Degrees of belief}\label{bayes}

\subsection{Boolean algebra models propositional logic}\label{logic}

Propositional logic studies finitary logical connectives like
``and'', ``or'', ``not'', which connect propositions.

Every proposition is either true or false.  Suppose some are
{\em known} to be true, some are {\em known} to be false, and
the truth values of some others are uncertain.  Call two propositions
$x$ and $y$ (conditionally) equivalent (given what is known)
if the proposition $[x\text{ if and only if }y]$ is known to be true.
It is easy to check that if $x_1$ is equivalent to $x_2$ and
$y_1$ is equivalent to $y_2$ then $[x_1\text{ and }y_1]$
is equivalent to $[x_2\text{ and }y_2]$,
$[x_1\text{ or }y_1]$ is equivalent to $[x_2\text{ or }y_2]$,
and $[\text{not }x_1]$ is equivalent to $[\text{not }x_2]$.
Therefore we can think of the three connectives
``and'', ``or'', ``not'' as acting on equivalence classes
rather than on propositions.

Any set of such equivalence classes of propositions that is
closed under these three connectives necessarily contains the
equivalence class, which we shall call 1, of propositions known
to be true, and the class, which we shall call 0, of propositions
known to be false.  If the truth values of an (equivalence
class of) proposition(s) is uncertain, then the set in question
also contains other classes than 0 and 1.  That set of equivalence
classes of propositions then constitutes a Boolean algebra with
the connectives ``and'', ``or'', and ``not'' in the roles of
meet, join, and complement.  The natural partial order of this
Boolean algebra makes $x\leq y$ precisely if the proposition
$[\text{if }x\text{ then }y]$ is known to be true.

\subsection{Intrinsic possibility versus epistemic possibility}

Possibility, like probability, can be either intrinsic or
epistemic.  To say it is possible that a card chosen randomly
from a deck will be an ace, could be taken to mean that
at least one ace is in the deck.
That is {\em intrinsic} possibility.  To say
it is possible that the card that was drawn yesterday was an
ace, could be taken to mean, not that some aces are in the deck,
but that it is not certain that none are.
That is {\em epistemic} possibility.
\begin{example}
Following the notation of \S\ref{logic}, we can say
that $x<y$ means it is possible that $y$ is true and $x$
is false, but it is necessary that $y$ is true if $x$ is true.
I was asked whether ``it is possible that $y$ is true and $x$
is false'' means
\begin{enumerate}
\item It is known that $y$ is possible without $x$; or
\item It is not known that $y$ is impossible without $x$.
\end{enumerate}
{\bf The punch line:} If possibility is regarded as intrinsic,
then (1) differs in meaning from (2), but if possibility is
regarded as epistemic, then there is no difference!
\end{example}
Henceforth we regard possibility as epistemic, not intrinsic.
That means, in particular, that we shall not speak of
$x$ as ``occurring'' or ``not occurring,'' but rather, as
we did earlier, of $x$ as being true or false, or as being
known to be true, known to be false, or uncertain.

Notice that the notation
$$
0\leq x<y\leq 1
$$
can be thought of as saying $x$ is closer to being
{\em known} to be false than $y$ is, or $y$ is closer
to being {\em known} to be true than $x$ is.  Consequently
we put a greater degree of belief in the truth of $y$
than in that of $x$.

\subsection{Some axioms of epistemic probability}

The last paragraph of the last section hints at an axiom
for epistemic probability theory: If $x$ is less (epistemically)
possible than $y$, then $x$ is less (epistemically) probable
than $y$.  In other words, for any assignment $P$ of probabilities
to propositions
\begin{equation}\label{firstaxiom}
\text{If }x<y\text{ then }P(x)<P(y).
\end{equation}
We shall also take it to be axiomatic that the less probable
$x$ is, the more probable $[\text{not }x]$ is, i.e.,
\begin{equation}\label{secondaxiom}
\text{If }P(x)<P(y)\text{ then }P(\text{not }x)>P(\text{not }y).
\end{equation}
A third axiom is very similar to the ``sure-thing principle''
stated by Leonard Jimmie Savage in \cite{savage}, pp.~21-2.
It says that if $x$ is no more probable than $y$ given that $z$
is true, and $x$ is no more probable than $y$ given that $z$
is false, then $x$ is no more probable than $y$ given no information
about whether $z$ is true or false.  In other words
\begin{equation}\label{thirdaxiom}
\text{If }P(x\mid z)\leq P(y\mid z)
\text{ and }P(x\mid \text{not }z)\leq P(y\mid\text{not }z)
\text{ then }P(x)\leq P(y)
\end{equation}
and ``$<$'' holds in the consequent if it holds in either of
the two antecedents.  (Savage's ``sure-thing principle'' spoke
of utilities rather than of probabilities.)  We do not understand
an expression like ``$P(\bullet\mid\bullet)$'' to mean anything
different from something like ``$P(\bullet)$''; we take all
probabilities to be conditional on some corpus of knowledge.
So in particular, (\ref{secondaxiom}) implies that
if $P(x\mid z)\leq P(y\mid z)$
then $P(\text{not }x\mid z)\geq P(\text{not }y\mid z)$.

It is unfortunate that, as things now stand, we must rely on
one more assumption about degrees of belief in uncertain
propositions -- that they are linearly ordered:
\begin{equation}\label{fourthaxiom}
\text{For all }x,y\text{ either }P(x)\leq P(y)
\text{ or }P(y)\leq P(x).
\end{equation}
This means we will have the conclusion we want for linearly
ordered scales and for scales that are Boolean algebras --
the two extreme cases -- but not for intermediate cases.

Clearly (\ref{firstaxiom}) says that assignments
of probabilities to propositions must satisfy
part 1(i) of Definition~\ref{BasicDefinitions}.
But (\ref{secondaxiom}) is weaker than
part 1(ii) of Definition~\ref{BasicDefinitions}.
If we can show that (\ref{secondaxiom}), (\ref{thirdaxiom}),
and (\ref{fourthaxiom}) require assignments of probabilities
to propositions to satisfy part 1(ii) of
Definition~\ref{BasicDefinitions} then we will know that
all such assigments must be basic scalings.  That is what
we do in the next section.

\subsection{Linearly ordered scales as probability assignments}

We want to show that (\ref{secondaxiom}), (\ref{thirdaxiom}),
and (\ref{fourthaxiom}) require assignments of probabilities
to propositions to satisfy part 1(ii) of
Definition~\ref{BasicDefinitions}.
Part 1(ii) of Definition~\ref{BasicDefinitions} speaks of
relative complementation.  In propositional logic, relative
complemenation is relative logical negation.
If $a\leq x\leq b$, meaning $a$ is a sufficient condition
for $x$, and $b$ is a necessary condition for $x$, then
the logical negation of $x$ relative to the interval
$[a,b]$ is the unique (up to logical equivalence) proposition
$u$ such that
\begin{enumerate}
\item $a$ is a sufficient condition for $u$, and
\item $b$ is a necessary condition for $u$, and
\item $u$ becomes equivalent to $[\text{not }x]$ once
it is learned that $b$ is true and $a$ is false.
\end{enumerate}
That proposition is $[a\text{ or }(b\text{ and not }x)]$,
or, equivalently (since $a$ logically entails $b$)
$[b\text{ and }(a\text{ or not }x)]$.  So the problem is
to show that (\ref{secondaxiom}), (\ref{thirdaxiom}), and
(\ref{fourthaxiom}) imply that if
$$
\begin{array}{ccccc}
a & \leq & x & \leq & b, \\
a & \leq & y & \leq & b,
\end{array}
$$
and $P(x)\leq P(y)$, then
$P(a\text{ or }[b\text{ and not }x])
\geq P(a\text{ or }[b\text{ and not }y])$.

If, to get a contradiction, we assume on the contrary that
$$
P(a\text{ or }[b\text{ and not }x])
\not\geq P(a\text{ or }[b\text{ and not }y])
$$
then
(\ref{fourthaxiom}) tells us that
$$
P(a\text{ or }[b\text{ and not }x])
<P(a\text{ or }[b\text{ and not }y]).
$$
The conjunction of this inequality with (\ref{thirdaxiom})
means we cannot have both
\begin{equation}\label{absurdinequality}
P(a\text{ or }[b\text{ and not }x]\mid b\text{ and not }a)
\geq P(a\text{ or }[b\text{ and not }y]\mid b\text{ and not }a)
\end{equation}
and
\begin{equation}\label{triviallyequal}
\begin{array}{cl}
& P(a\text{ or }[b\text{ and not }x]\mid
\text{\underline{not} }\left\{b\text{ and not }a\right\}) \\  \\
= & P(a\text{ or }[b\text{ and not }y]\mid
\text{\underline{not} }\left\{b\text{ and not }a\right\}).
\end{array}
\end{equation}
The equality in (\ref{triviallyequal}) is trivially true
because the condition $[\text{not }(b\text{ and not }a)]$
renders impossible the propositions whose probability is
being taken.  Therefore (\ref{absurdinequality}) must be false.
So, by (\ref{fourthaxiom}), we must have
\begin{equation}\label{toocomplicated}
P(a\text{ or }[b\text{ and not }x]\mid b\text{ and not }a)
<P(a\text{ or }[b\text{ and not }y]\mid b\text{ and not }a).
\end{equation}
Given $[b\text{ and not }a]$, the propositions
$[a\text{ or }(b\text{ and not }x)]$ and
$[a\text{ or }(b\text{ and not }y)]$ simplify to
$[\text{not }x]$ and $[\text{not }y]$ respectively,
and (\ref{toocomplicated}) simplifies to
\begin{equation}\label{simplified}
P(\text{not }x\mid b\text{ and not }a)
<P(\text{not }y\mid b\text{ and not }a).
\end{equation}
We must also have
\begin{equation}\label{triviallyequal2}
P(\text{not }x\mid
\text{\underline{not} }\left\{b\text{ and not }a\right\})
=P(\text{not }y\mid
\text{\underline{not} }\left\{b\text{ and not }a\right\})
\end{equation}
because, given the condition $[\text{not }(b\text{ and not }a)]$,
the two propositions $[\text{not }x]$ and $[\text{not }y]$ are
equivalent to each other.  Conjoining
(\ref{simplified}), (\ref{triviallyequal2}), and (\ref{thirdaxiom}), 
we conclude that
$P(\text{not }x)<P(\text{not }y)$.
In view of (\ref{secondaxiom}) and our assumption
that $P(x)\leq P(y)$, this is impossible.

We conclude that assignments of linearly ordered probabilities
to uncertain propositions should be scalings.  As scalings, they
must satisfy all of our results on addition, dual-addition,
subtraction, relative complementation, modularity, and de-Morganism.

\subsection{Finiteness of information-content in propositions}

When would a scaling need to be continuous in order to model
properly the phenomenon of assignment of degrees of belief to
uncertain propositions?

Suppose a subset $\mathbb{X}$ of some Boolean algebra
$\mathbb{A}$ of propositions is closed under ``and''
(i.e., $[x\text{ and }y]\in\mathbb{X}$
for any $x,y\in\mathbb{X}$) and satisfies
$\bigwedge\mathbb{X}=0$ and $x>0$ for every $x\in\mathbb{X}$.
The simplest example is the Boolean algebra $\mathbb{A}$
that is freely generated by $x_1,x_2,x_3,\dots$, i.e.,
the set of all propositions constructed from $x_1,x_2,x_3,\dots$
by using only finitely many occurences of ``and'', ``or'',
and ``not'', and $\mathbb{X}$ is the set
$\set{x_1,x_2,x_3,\dots}$ of generators.
If a probability $\alpha$ satisfies $0<\alpha\leq P(x)$
for every $x\in\mathbb{X}$, then it would seem appropriate
to consider $\alpha$ to be a probability assigned only to
propositions that convey an amount of information that is
infinite by comparison to that conveyed by any $x\in\mathbb{X}$.
Closure of $\mathbb{X}$ under ``and'' is the same as closure
of a family of clopen subsets of the Stone space under finite
intersections.  Consequently, compactness of the Stone space
implies that for some 2-valued homomorphism $\varphi$ we have
$\varphi(x)=1$ for every $x\in\mathbb{X}$.  The probability
$\alpha$ described above would then be a counterexample to
the continuity of $P$ at $\varphi$.  So exclusion from a scale,
of probabilities assigned only to propositions that convey
an infinite amount of information, amounts to continuity of
the assignment of probabilities at every 2-valued homomorphism.
In view of Theorems~\ref{archimedean} and~\ref{measurable}, then,
an insistence on finite information content in propositions,
implies that the scale on which the probabilities are measured
is measurable.  In other words, we may take the probabilities
to be real numbers, and the scale to be $[0,1]\subseteq\mathbb{R}$.

\subsection{Multiplication}\label{productrule}

From (\ref{thirdaxiom}) it follows that
\begin{equation}\label{simpleinequality}
P(x\text{ and }z)\leq P(y\text{ and }z)
\text{ if and only if }P(x\mid z)\leq P(y\mid z).
\end{equation}
By symmetry
the same is true if ``$\geq$'' replaces both
occurences of ``$\leq$'', and consequently also
of ``$=$'' replaces both.

If, as in \S\ref{secmultiplication}, we assume linear ordering,
dividedness, and Archimdeanism, then there exists a maximum
non-negative integer $n_1$ such that
\begin{equation}\label{dividingprobabilities}
P(z)\underbrace{
-P(x\text{ and }z)-\cdots\cdots-P(x\text{ and }z)}_{n_1}
\end{equation}
exists, and the difference is less than $P(x\text{ and }z)$.
Because of (\ref{simpleinequality}), this integer $n_1$
must be the same as the smallest non-negative integer $n$
such that\begin{equation}\label{dividingconditionalprobabilities}
1\underbrace{
-P(x\mid z)-\cdots\cdots-P(x\mid z)}_{n}
\end{equation}
exists.  This latter difference must be less than $P(x\mid z)$.
The process can be iterated according to the algorithm described
in \S\ref{secmultiplication}, and (\ref{simpleinequality}) tells
us at each step that the entries $n_1,n_2,n_3\dots$ in the
formal continued fraction expansion of the formal quotient
$P(x\text{ and }z)/P(z)$ are the same as those in the
formal continued fraction expansion of the formal quotient
$P(x\mid z)/1$.  Thus we have:
\begin{theorem}
If probabilities are measured on a scale that is
linearly ordered, divided, and Archimedean, then for
any propositions $x$ and $z$,
$$
P(x\text{ and }z)=P(x\mid z)P(z).
$$
\end{theorem}

\begin{acknowledgments}
James M.~Dickey, Karel Prikry, and Gian-Carlo Rota all
encouraged me in different ways to work on this topic.
From a usenet posting by Wlodzimierz Holsztynski I first
heard of Kleene algebras, and he assisted in my understanding
them in a subsequent exchange of e-mail.  Richard Stanley's
comment on the confusing nature of a bit of nomenclature
inspired the term ``additive relative complement.''
Sara Billey suggested a need to clarify the ``two-page crash
course on nonstandard analysis'' that appears in
Example~\ref{fullnonarch}.  Timothy Chow and Daniel
Lueking independently suggested the term ``stably Archimedean.''
\end{acknowledgments}

\end{article}
\end{document}